\documentclass[english, 12pt,leqno]{amsart}

\usepackage{amsmath}
\usepackage{amsfonts}
\usepackage{amssymb}

\title[Max Dehn]
{Topologie,  th\'eorie des groupes 
\\
et probl\`emes de d\'ecision
\\
---
\\
c\'el\'ebration d'un article de Max Dehn 
\\
de 1910
}

\author{Pierre de la Harpe
}

\thanks{L'auteur a tr\`es largement profit\'e 
de nombreuses conversations avec Claude Weber, 
au si\`ecle pass\'e
pour son \'education topologique en g\'en\'eral,
et plus r\'ecemment
pour la r\'edaction de ce texte en particulier.
Pour plusieurs remarques utiles sur des versions pr\'eliminaires du texte
et des encouragements \`a divers stades du travail,  
je remercie aussi 
Michel Boileau, Martin Bridson, Etienne Ghys,
Cameron Gordon, Tatiana Nagnibeda, Mark Powell,
Jean-Philippe Pr\'eaux et Thierry Vust.}

\address{Pierre de la Harpe~:
Section de math\'ematiques, 
Universit\'e de Gen\`eve, 
C.P.~64, 
CH--1211 Gen\`eve 4. 
}
\email{Pierre.delaHarpe@unige.ch}


\begin{document}

\begin{abstract} 
This paper, in French, is a celebration of Max Dehn,
and an essay of describing some of his results
published in the beginning of the 1910's, and their offspring.
It has been written up for a winter school in Les Diablerets,
March 7--12, 2010:

\emph{
Geometry, topology and computation in groups,
}\par\emph{
100 years since Dehn's Decision Problems.
}

\vskip.2cm\noindent
\textsc{R\'esum\'e.}
Ce texte est une c\'el\'ebration de Max Dehn
et un essai de mise en perspective
de quelques-uns de ses r\'esultats publi\'es au d\'ebut des ann\'ees 1910.
Il a \'et\'e r\'edig\'e \`a l'occasion d'une \'ecole d'hiver aux Diablerets
du 7 au 12 mars 2010~:

\emph{
Geometry, topology and computation in groups,
}\par\emph{
100 years since Dehn's Decision Problems. 
}  

\end{abstract}

\maketitle

\begin{center}
A para\^{\i}tre dans la \emph{Gazette des math\'ematiciens}
\end{center}

\section{El\'ements biographiques}
\label{section1}

Max Dehn est n\'e \`a Hambourg en 1878, 
et mort \`a Black Mountain (Caroline du Nord, U.S.A.) en 1952.
Il fut \'etudiant de Hilbert \`a G\"ottingen en 1899
et acheva en 1900 une th\`ese sur 
les fondements de la g\'eom\'etrie.
\par

En 1900 \'egalement, il r\'esolut  le \emph{troisi\`eme probl\`eme de Hilbert},
en montrant que deux t\'etra\`edres de m\^eme volume dans $\mathbf R^3$
ne sont pas n\'ecessairement \'equid\'ecomposables.
Il en r\'esulte que, contrairement \`a 
la th\'eorie des aires des polygones dans $\mathbf R^2$,
la th\'eorie  des volumes des poly\`edres dans $\mathbf R^3$
doit reposer
sur la notion de limite, ou sur son anc\^etre
qui est la m\'ethode d'exhaustion d'Eudoxe (-408 -- -355),
Euclide ($\sim$ -325 -- $\sim$ -265)
et Archim\`ede ($\sim$ -287 -- -212).
[Bien s\^ur, le sujet n'est pas clos~! \cite{Zeem--02}.]

Apr\`es un bref s\'ejour comme assistant \`a Karlsruhe,
Dehn fut privat-docent \`a M\"unster jusqu'en 1911~;
durant cette p\'eriode, sous l'influence de Poul Heegaard et Henri Poincar\'e,
il commen\c ca \`a s'int\'e\-resser \`a la topologie et \`a la th\'eorie des groupes.
En 1907, Dehn et Heegaard publi\`erent le panorama
de l'\emph{Enzyklop\"adie der mathematischen Wissenschaften}
sur l'Analysis Situs \cite{DeHe--07}~; 
ce texte contient une  
d\'efinition et une classification des surfaces
du point de vue de la topologie combinatoire.
Et c'est en 1910 que Dehn publia le premier des articles majeurs
o\`u apparaissent les \emph{probl\`emes de d\'ecision}.

Dehn fut ensuite professeur assistant \`a Kiel  (1911--13), 
professeur ordinaire \`a Breslau (1913--21),
avec une interruption due au  service arm\'e (1915--18),
et successeur de Bieberbach \`a Francfort (1922--35).
Le s\'eminaire d'histoire des math\'ematiques qu'il y dirigea
semble avoir \'et\'e un grand moment pour tous les participants
\cite{Sieg--64}. 
Il fut patron de th\`ese de 
\begin{itemize} 
\item[-]
Hugo Gieseking (1912), 
qui construisit une vari\'et\'e hyperbolique de dimension $3$,
non orientable, dont Thurston nota que le rev\^etement orientable
est isom\'etrique au compl\'ement d'un n\oe{}ud de huit, 
avec sa structure de vari\'et\'e hyperbolique 
(re?)d\'ecouverte en 1975 par Robert Riley (voir par exemple \cite{Miln--82})~;
\item[-]
Jakob Nielsen (1913), l'un des fondateurs de la th\'eorie combinatoire des groupes,
dont les travaux sur les diff\'eomorphismes de surfaces sont de grande importance
\cite{Thur--76}~;
\item[-]
Ott-Heinrich Keller (1929),
dont la th\`ese sur les pavages de $\mathbf R^n$ par des cubes
contient une conjecture aux nombreuses ramifications,
autant que je sache toujours ouverte lorsque $n=7$
(voir le chapitre 7 de \cite{Zong--05})~;
Keller formula aussi la ``conjecture du jacobien'', ch\`ere \`a Abhyankar~;
\item[-]
Wilhelm Magnus (1931), qui d\'emontra dans sa th\`ese le
``Freiheitssatz'' formul\'e par Dehn, sur les groupes \`a un relateur, 
et qui s'illustra plus tard dans de multiples domaines \cite{Magn--94}~;
\item[-]
Ruth Moufang (1931), dont la th\`ese portait sur la g\'eom\'etrie projective
et \`a qui on doit des contributions majeures sur
les structures alg\'ebriques non associatives.
Pour illustration~: dans sa construction du ``monstre''
(le plus grand groupe fini simple sporadique),
Conway a utilis\'e une ``Moufang loop'' construite par Parker~;
\end{itemize}
liste \`a laquelle le
``Mathematics Genealogy Project'' ajoute
Herbert Fuss (1913), Wilhelm Schwan (1923), Max Frommer (1928), et Joseph Engel (1949).
\par

En \'et\'e 1935, Dehn fut d\'emis de son poste pour raison d'ascendance juive.
Il continua \`a (sur)vivre en Allemagne jusqu'en 1939, 
avant de s'\'echap\-per vers les Etats-Unis
via Copenhague, Trondheim, la Finlande, le Transsib\'erien et le Japon.
Apr\`es divers postes de courtes dur\'ees,
il termina sa carri\`ere  au petit coll\`ege de Black Mountain, 
dans l'est des Etats-Unis (1945--52).

Il y fut l'unique professeur de math\'ematiques~;
il y enseigna aussi la philosophie, le latin et le grec. 
Ce coll\`ege, fond\'e en 1933, 
avait des objectifs p\'edagogiques originaux 
et ambitieux\footnote{``The college sought 
to educate the whole student -- head, heart and hand -- through studies, 
the experience of living in a small community and manual work.''
Ceci est bien plus sur le site
{\tt http://www.bmcproject.org/index.htm}}
qui plurent \`a Dehn.
Son salaire mensuel initial  \'etait de 40 \$, plus log\'e, nourri, blanchi.
(C'\'etait mieux que l'offre initiale, qui \'etait de 20 \$~;
mais Siegel raconte qu'il y eut des p\'eriodes pendant lesquelles
l'argent de poche, en compl\'ement du g\^{\i}te et du couvert,
\'etait limit\'e \`a 5 \$ par mois.)
Le coll\`ege avait des probl\`emes financiers et ferma en 1956.
Dehn y avait eu de longues conversations avec un coll\`egue architecte,
aujourd'hui c\'el\`ebre pour ses constructions utilisant des poly\`edres r\'eguliers~: Richard Buckminster Fuller \cite{BuFu}.
Dehn \'etait aussi naturaliste amateur et randonneur enthousiaste  \cite{Sher--94}.

Pour en lire davantage sur la vie et la production math\'ematique de Dehn, 
voir \cite{Magn--78}, \cite{Stil--99} et \cite{Daws--02}.

\section{Les trois probl\`emes de d\'ecision de Dehn}
\label{section2}

Dans trois articles publi\'es de 1910 \`a 1912,
Max Dehn a formul\'e et \'etudi\'e les trois \emph{probl\`emes de d\'ecision} suivants
qui sont fondamentaux en th\'eorie combinatoire des groupes.
\begin{itemize}
\item[(WP)] \textbf{Le probl\`eme du mot}~:
pour un groupe $G$ engendr\'e par des \'el\'ements
$s_1, s_2, \hdots, s_n$,
\emph{ trouver une m\'ethode qui permette de d\'ecider
en un nombre fini de pas 
si deux produits des op\'erations $s_i$ de $G$ sont \'egaux,
en particulier si un tel produit d'op\'erations est \'egal \`a l'identit\'e.}

\item[(CP)] \textbf{Le probl\`eme de conjugaison}~: 
pour $G$ et les $s_i$ comme ci-dessus,
\emph{trouver une m\'ethode qui permette de d\'ecider
en un nombre fini de pas,
\'etant donn\'e deux substitutions $S$ et $T$ de $G$,
s'il existe une troisi\`eme substitution $U$ de $G$
telle que $S = UTU^{-1}$,
c'est-\`a-dire si $S$ est un conjugu\'e de $T$.}

Ces deux probl\`emes sont formul\'es dans \cite{Dehn--10}, 
ici traduits de la version anglaise, page 95 de \cite{DeSt--87}.
La formulation de \cite{Dehn--11} est l\'eg\`erement diff\'erente.

\item[(IsP)] \textbf{Le probl\`eme d'isomorphisme}~: 
\emph{\'etant donn\'e deux groupes, d\'ecider s'ils sont isomorphes ou non
(et de plus si une correspondance donn\'ee entre
les g\'en\'erateurs d'un groupe et les \'el\'ements de l'autre
est ou n'est pas un isomorphisme).}

Probl\`eme formul\'e dans \cite{Dehn--11}, 
page 134 de \cite{DeSt--87}.

\end{itemize}
Dans une formulation concise, reprise de \cite{Mill--92}~:
\begin{itemize}
\item[-] 
$\operatorname{WP}(G) \, = \, (?w \in G)(w =_G 1) ,$
\item[-]
$\operatorname{CP}(G) \, = \, (?u,v \in G)(\exists x \in G)(x^{-1}ux =_G v) ,$
\item[-] 
$\operatorname{IsP} \,  = \, 
(?\pi_1, \pi_2 \hskip.2cm \text{pr\'esentations finies})(gp(\pi_1) \approx gp(\pi_2)).$
\end{itemize}
Ici,  ``?...'' est une esp\`ece de quantificateur, signifiant par exemple pour $WP(G)$,
``le probl\`eme de d\'ecider, pour un $w \in G$ arbitraire, si oui ou non $w =_G 1$''.
Les deux premiers probl\`emes de d\'ecision de Dehn 
ont des formulations topologiques, 
et le troisi\`eme est \'egalement li\'e \`a une question topologique~:
\begin{itemize}
\item[(WP$_{\text{top}}$)] 
dans un espace topologique de groupe fondamental $G$,
est-ce qu'un lacet ``donn\'e'' est contractile sur un point~?
\item[(CP$_{\text{top}}$)] 
est-ce que deux lacets ``donn\'es''
sont librement homotopes~?
\item[(IsP$_{\text{top}}$)] 
est-ce que deux espaces ``donn\'es'' 
sont homotopiquement \'equi\-valents~?
\end{itemize}
(les guillemets parce que nous n'allons pas essayer
de pr\'eciser ce que ``donn\'e''  veut dire).

\medskip

Quelques premi\`eres remarques~:

(i) A priori, la r\'eponse au probl\`eme du mot
d\'epend de la pr\'esentation consid\'er\'ee du groupe.
Toutefois, si
$\langle S_1 \hskip.1cm \vert \hskip.1cm R_1 \rangle$ et 
$\langle S_2 \hskip.1cm \vert \hskip.1cm R_2 \rangle$
sont deux pr\'esentations finies d'un m\^eme groupe $G$, 
un argument simple montre \emph{qu'il existe}
``une m\'ethode qui permette de d\'ecider  ...''
pour une pr\'esentation 
si et seulement s'il en existe une pour l'autre.
Il peut n\'eanmoins \^etre plus d\'elicat \emph{d'exhiber} une m\'ethode pour
$\langle S_2 \hskip.1cm \vert \hskip.1cm R_2 \rangle$
\`a partir d'une m\'ethode pour 
$\langle S_1 \hskip.1cm \vert \hskip.1cm R_1 \rangle$.
\par
Consid\'erons par exemple un groupe fini $G$,
et sa table de multiplication vue comme une pr\'esentation
$\langle G \hskip.1cm \vert \hskip.1cm R_{\text{mult}} \rangle$,
o\`u $ R_{\text{mult}}$ d\'esigne l'ensemble des relations
$abc^{-1} = 1$, avec $a,b \in G$ et $c$ le produit $ab$~;
il est tout \`a fait banal d'\'ecrire un algorithme r\'esolvant le probl\`eme du mot
pour cette pr\'esentation.
Cela l'est en revanche beaucoup moins pour une pr\'esentation arbitraire
du m\^eme groupe fini $G$.
L'une des m\'ethodes utilisables remonte \`a un article
de Todd et Coxeter \cite{ToCo--36}~; 
voir aussi  \cite{John--80}, \cite{Cann--02} et \cite{HoEO--05}.
\par
Des remarques analogues valent pour le probl\`eme de conjugaison.

(ii) Si le probl\`eme de conjugaison est r\'esoluble pour un groupe $G$,
alors le probl\`eme du mot l'est aussi,
puisqu'un \'el\'ement de $G$ est \'egal \`a l'identit\'e
si et seulement s'il est conjugu\'e \`a l'identit\'e.
La r\'eciproque n'est pas vraie.

(iii) Voici un exemple \`a propos du probl\`eme d'isomorphisme~:
le groupe donn\'e par la pr\'esentation
\begin{equation*}
\langle s,t \hskip.1cm \vert \hskip.1cm s^3t, t^3, s^4 \rangle
\end{equation*}
est le groupe \`a un \'el\'ement.
Ce n'est peut-\^etre pas tout \`a fait \'evident, mais c'est vrai~:
$s = ss^3t = s^4t = t$ et $1 = s^3tt^{-3} = s = t$.

(iv) Les articles de Dehn sont bien ant\'erieurs 
aux d\'efinitions pr\'ecises de mots comme ``algorithmes''
ou ``proc\'ed\'es'' ou ``m\'ethodes''~;
en fait, c'est un des grands succ\`es de la logique math\'ematique
que d'avoir donn\'e un sens pr\'ecis aux probl\`emes de Dehn
(Post, Church, Turing, G\"odel, ... dans les ann\'ees 1930).

Ce n'est que dans les ann\'ees 1950
qu'il fut d\'emontr\'e que ces probl\`emes sont insolubles
pour certains groupes de pr\'esentation finie~:
Piotr Sergue\"{\i}evitch Novikov\footnote{A ne pas confondre avec
son fils Sergue\"{\i} Petrovitch Novikov,
laur\'eat de la m\'edaille Fields en 1970,
pour des travaux portant entre autres sur le cobordisme,
les feuilletages, et l'invariance topologique des classes de Pontrjagin
d'une vari\'et\'e diff\'erentiable.
} 
d\'emontra que le probl\`eme de conjugaison est insoluble (1954),
Boone et Novikov que le probl\`eme du mot est insoluble (1955),
et Adjan et Rabin que le probl\`eme d'isomorphisme est insoluble (1958).
Fridman montra que le probl\`eme du mot peut \^etre r\'esoluble
et le probl\`eme de conjugaison insoluble pour le m\^eme groupe (1960)~;
Miller (1971) a montr\'e que c'est m\^eme le cas pour un groupe
qui s'ins\`ere dans une suite exacte courte
\`a noyau et quotient libres (th\'eor\`eme 4.8 de \cite{Mill--92}).
Pour en savoir plus, voir  \cite{Mill--71}, \cite{Mill--92} et \cite{Stil--82}.

(v) Les probl\`emes de d\'ecision furent formul\'es par Dehn
dans le contexte de la \emph{topologie de dimension $3$}
et de la \emph{th\'eorie des n\oe{}uds}\footnote{D\`es 1867,
le physicien \'ecossais Peter Guthrie Tait avait 
commenc\'e \`a dresser des tables de n\oe{}uds~;
mais il n'y a pas de topologie dans les travaux de Tait,
et a fortiori pas de groupes de n\oe{}uds.
La th\'eorie des n\oe{}uds  intervient aussi dans l'\'etude des
singularit\'es isol\'ees des fonctions de deux variables complexes,
comme cela est apparu vers 1905 dans les recherches de Wirtinger~;
celles-ci n'eurent qu'une diffusion limit\'ee \`a l'\'epoque.
Ce n'est que bien plus tard que la th\'eorie des n\oe{}uds
devint  ``respectable'' pour cette raison,
lorsque les recherches de Wirtinger furent reprises par
Brauner (dans son habilitation de 1928 sous la supervision de Wirtinger),
K\"ahler, Zariski et Burau~; 
voir \cite{Eppl--95}, d\'ej\`a cit\'e,
ou les pages 318-320 de \cite{Eppl--99a}.
Ce lien entre n\oe{}uds et singularit\'es ne semble pas avoir jou\'e
de r\^ole dans les motivations de Dehn.}.
Il faut garder en t\^ete que la topologie alg\'ebrique
\'etait de cr\'eation r\'ecente~:
on peut la dater de 1895, 
ann\'ee de parution de l'\emph{Analysis situs} 
\cite{Poin--95}\footnote{Plus 
une \emph{Note} pr\'eliminaire de 1892 
et les cinq \emph{compl\'ements}
publi\'es entre 1899 et 1904 \cite{Poin--04},
notamment suite \`a la th\`ese de Heegaard (1898)
qui mettait en \'evidence une erreur de Poincar\'e
(oubli de la partie de torsion des groupes d'homologie).
On peut souligner la rapidit\'e d'une certaine \'evolution~:
les motivations explicites de Poincar\'e \'etaient li\'ees \`a l'analyse~:
courbes d\'efinies par des \'equations diff\'erentielles,
probl\`eme des trois corps, fonctions multi-valu\'ees \`a deux variables,
p\'eriodes des int\'egrales multiples, calculs de perturbation, ...~;
mais la topologie est tr\`es vite devenue une discipline autonome
de l'analyse, c'est flagrant dans les articles de Dehn.
}.
\par

L'int\'er\^et de Dehn pour la topologie
avait peut-\^etre \'et\'e \'eveill\'e (entre autres) par
sa r\'esolution d'un ``exercice'' que lui avait donn\'e Hilbert,
c'est-\`a-dire par sa d\'emonstration du th\'eor\`eme de Jordan 
pour les polygones plans.
Ce travail de Dehn ne fut pas publi\'e \cite{Gugg--77}. 
\par

Dehn avait esp\'er\'e
r\'esoudre la \emph{conjecture de Poincar\'e},
concernant les vari\'et\'es closes de dimension~$3$
\`a groupe fondamental trivial. 
Plus pr\'ecis\'ement, le 12 f\'evrier 1908, Dehn avait envoy\'e \`a Hilbert
un article contenant un ``r\'esultat'' sur $\mathbf R^3$ 
\'equivalent \`a la conjecture de Poincar\'e, 
article soumis pour publication dans les ``G\"ottinger Nachrichten''.
Mais Tietze signala une erreur dans l'argument de Dehn,
qui retira son article par une lettre du 16 avril de la m\^eme ann\'ee.
(D'apr\`es la page 388 de \cite{Eppl--95}.)
Ainsi Dehn fut-il ``la premi\`ere victime de la conjecture de Poincar\'e''
(selon une formule de K. Volkert, 1996).
\par

(vi) Les lecteurs au courant des articles de l'\'epoque notent que
le probl\`eme d'isomorphisme avait d\'ej\`a \'et\'e formul\'e par Tietze en 1908,
mais sans que Tietze lui donne l'importance que Dehn donna \`a ses formulations.
\par

Bien que Tietze et Dehn comptent tout deux parmi les fondateurs de la topologie, 
et en particulier de la th\'eorie des n\oe{}uds qui donna \`a Dehn l'occasion
de formuler ses probl\`emes de d\'ecision,
les influences qu'ils eurent chacun sur l'autre \`a cette \'epoque
semblent avoir \'et\'e limit\'ees 
(nonobstant ce qui est rapport\'e plus haut).
Tietze, \'etudiant \`a Vienne de 1898 jusqu'\`a son habilitation en 1908
(avec une ann\'ee \`a Munich), \'etait sous l'influence de Wirtinger,
et par lui dans le sillage de Klein,
alors que Dehn, \'etudiant \`a G\"ottingen,
\'etait marqu\'e par l'influence directe et bien diff\'erente de Hilbert.
Voir \cite{Eppl--95}, en particulier la page 395.
\par

(vii) En th\'eorie des groupes, la formulation  des probl\`emes de d\'ecision,
due \`a Dehn,  fut un moment cl\'e.
Mais ce ne sont pas les seuls probl\`emes de d\'ecision qui puissent \^etre formul\'es~!
Ne mentionnons ici que le dixi\`eme probl\`eme de Hilbert~:
\emph{De la possibilit\'e de r\'esoudre une \'equation diophantienne. On donne une \'equation de Diophante \`a un nombre quelconque d'inconnues et \`a coefficients entiers rationnels : on demande de trouver une m\'ethode par laquelle, au moyen d'un nombre fini d'op\'erations, on pourra distinguer si l'\'equation est r\'esoluble en nombres entiers rationnels.}
La r\'eponse, n\'egative, r\'esulte des travaux de
Martin Davis, Yuri Matiyasevich, Hilary Putnam et Julia Robinson~;
le coup de gr\^ace fut donn\'e par Matiyasevich en 1970. 

\medskip
\begin{center}
* * * * * * * * * * * * *
\end{center}

Toujours est-il que, peu avant 1910,
Dehn ma\^{\i}trisait aussi bien que quiconque \`a l'\'epoque
des sujets que nous \'enum\'erons comme suit~;
ceci  dans le vocabulaire d'aujourd'hui,
donc en termes parfaitement anachroniques.

\medskip

($\spadesuit$) 
La notion de \emph{pr\'esentation d'un groupe}
\begin{equation*}
G \, = \, \langle S \hskip.1cm \vert \hskip.1cm R \rangle \, = \, 
 \langle s_1, \hdots \hskip.1cm \vert \hskip.1cm r_1, \hdots \rangle 
\end{equation*}
donn\'e comme quotient $G = F/R$ 
d'un groupe libre $F$ de base $S = \{s_1, \hdots\}$
par le sous-groupe normal engendr\'e par les ``relations'' $R = \{r_1, \hdots \}$.
La notion fut pr\'ecis\'ee par Walther von Dyck (1882), 
suite \`a des travaux sur les groupes discontinus
apparaissant en th\'eorie des fonctions de variables complexes,
dont ceux de son ma\^{\i}tre Felix Klein.

\medskip

($\heartsuit$)
Le \emph{groupe fondamental} d'un espace topologique,
qui fut introduit par Poincar\'e en 1895 dans son 
\emph{Analysis situs}\footnote{Chez 
Poincar\'e, la premi\`ere d\'efinition \emph{explicite} du groupe fondamental
est en termes de rev\^etements [Note de 1892], 
et la seconde en termes de classes d'homotopie de lacets \cite{Poin--95}~;
de plus, on peut voir le groupe fondamental entre les lignes de \cite{Poin--83},
au sujet de l'uniformisation des surfaces de Riemann
(merci \`a Etienne Ghys pour cette observation).
Quelques d\'etails de plus aux pages 374-376 de \cite{Eppl--95}.
Pour lire Dehn, il faut adopter le point de vue des classes de lacets.
}.
En particulier, le groupe fondamental d'une vari\'et\'e compacte
est un groupe de pr\'esentation finie, comme cela appara\^{\i}t
clairement dans l'habilitation de Tietze (1908)~;
il en est de m\^eme pour le groupe fondamental
d'un complexe simplicial fini 
(et m\^eme d'un CW-complexe fini).

\medskip

($\diamondsuit$)
En particulier, un 
\emph{n\oe{}ud, c'est-\`a-dire une courbe ferm\'ee simple $K$ 
dans l'espace usuel,}
d\'etermine un groupe $G_K$,
le groupe fondamental du compl\'ementaire du n\oe{}ud.
\par

On peut facilement \'ecrire une  pr\'esentation finie de $G$
en termes d'un diagramme plan $\mathcal D$ de $K$.
On conna\^{\i}t ainsi la
\emph{pr\'esentation de Dehn}\footnote{Si j'en crois ce qu'affirment les lecteurs 
plus perspicaces que moi,
cette \emph{pr\'esentation de Dehn} est expliqu\'ee
au \S~3 du chapitre II de \cite{Dehn--10}.},
dont les g\'en\'erateurs sont en bijection avec
les composantes connexes born\'ees du compl\'ementaire
de $\mathcal D$ dans le plan 
et les relations en bijection avec les points doubles de $\mathcal D$
(voir par exemple l'appendice de \cite{Kauf--83}),
et la \emph{pr\'esentation de Wirtinger}, 
qui remonte \`a un expos\'e de  Wirtinger de 1905,
et qu'on trouve dans presque tous les expos\'es
de th\'eorie des n\oe{}uds (par exemple au \S~5 de \cite{Rham--69}
ou au chapitre~3 de \cite{BuZi--85}).
M\^eme si l'expos\'e de Wirtinger est bien ant\'erieur  \`a l'article de Dehn,
``apparemment, Dehn n'\'etait pas conscient du r\'esultat de Wirtinger''
(citation de \cite{Magn--78}). 
\par

A propos du probl\`eme d'isomorphisme,
noter que deux diagrammes $D$, $D'$ d'un m\^eme n\oe{}ud $K$
donnent en g\'en\'eral lieu \`a des pr\'esentations distinctes 
$\pi_D$, $\pi_{D'}$ du m\^eme groupe $G_K$,
et il peut \^etre non banal d'\'ecrire un argument combinatoire
montrant qu'on peut transformer une pr\'esentation en l'autre.
Plus g\'en\'eralement, 
il peut \^etre bien difficile de distinguer ``\`a l'oeil nu''
si deux diagrammes $D$, $D'$ d\'efinissent ou non le m\^eme n\oe{}ud,
ou si les pr\'esentations correspondantes d\'efinissent ou non 
des groupes isomorphes.
Par exemple, Tietze s'est amus\'e \`a dessiner 
deux diagrammes de n\oe{}uds $D$ et $D'$
\`a $48$ (sauf erreur) croisements chacun,
qui ont l'air identiques
(ils se distinguent en un seul croisement),  
mais $D$ repr\'esente le n\oe{}ud de tr\`efle et $D'$ le n\oe{}ud trivial
\cite{Tiet--42}.
J'imagine qu'il ne serait pas imm\'ediat de d\'ecider par force brutale
que les pr\'esentations correspondantes, 
\`a $48$ g\'en\'erateurs et $48$ relations,
d\'efinissent deux groupes non isomorphes.

\medskip

($\clubsuit$) 
Vers 1909-1910, 
Dehn avait donn\'e et r\'edig\'e au moins deux chapitres de cours
(non publi\'es de son vivant 
mais r\'ecemment parus en traduction anglaise dans \cite{DeSt--87}),
sur la th\'eorie des groupes \cite{Dehn--a}
et la topologie des surfaces  \cite{Dehn--b}.
Il en ressort entre autres que Dehn,
comme d'ailleurs beaucoup de ses contemporains, 
\'etait tout \`a fait \`a l'aise
avec la \emph{g\'eom\'etrie hyperbolique} (plane).
En particulier, l'\'etude des ``groupes fuchsiens''
conduit \`a des pavages du plan hyperbolique
par des polygones fondamentaux,
pavages dont les graphes duaux sont pr\'ecis\'ement
des ``Dehn's Gruppenbilder'' 
(ou ``graphes de Cayley'', voir ci-dessous le chapitre \ref{section4}).
Dans ce sens, les ``Dehn's Gruppenbilder''
doivent beaucoup plus \`a Dyck, Klein et Fricke qu'\`a Cayley.
Voir les figures de \cite{KlFr--90} et \cite{FrKl--97}~;
voir aussi les figures des chapitres XVIII et XIX de \cite{Burn--11}.

\medskip
\begin{center}
* * * * * * * * * * * * *
\end{center}

Le premier des articles majeurs de Dehn qui nous int\'eressent ici
\cite{Dehn--10}
parut dans le premier cahier du volume 69 des \emph{Mathematische Annalen},
``Ausgegeben am 23.~Juni 1910''.
On y trouve  les nouveaut\'es suivantes~;
dans la liste ci-dessous, (n) se r\'ef\`ere au chapitre n du pr\'esent texte, dont
les chapitres n$^{bis}$ sont consacr\'es \`a certains d\'eveloppements
plus r\'ecents.
\begin{itemize}
\item[(2)]
l'\'enonc\'e du probl\`eme du mot et du probl\`eme de conjugaison~;
\item[(3)]
le lemme de Dehn, et le crit\`ere de trivialit\'e qui en d\'ecoule pour les n\oe{}uds~;
\item[(4)]
un diagramme de Cayley (ou plut\^ot de Dehn~!) du groupe du n\oe{}ud de tr\`efle~;
\item[(5)]
une construction de $3$-sph\`eres d'homologie, 
dont l'une a un groupe fondamental fini non r\'eduit \`a un \'el\'ement~;
\item[(7)]
tout groupe de pr\'esentation finie est le groupe fondamental d'un $2$-complexe
(et aussi d'une $4$-vari\'et\'e close, voir le chapitre III de \cite{Dehn--11}) ~;
\end{itemize}
excusez du peu.
A peine plus tard, Dehn publiera
\begin{itemize}
\item[(6)]
des algorithmes pour r\'esoudre le probl\`eme du mot 
et le probl\`eme de conjugaison dans les groupes fondamentaux 
des surfaces orientables~; 
groupes sur lesquels Dehn consid\`ere la \emph{m\'etrique des mots}
et note qu'elle peut remplacer avantageusement la m\'etrique hyperbolique~;
\item[(-)]
une d\'emonstration du fait que les n\oe{}uds de tr\`efle gauche et droit
ne sont pas isotopes (exposition de cette d\'emonstration 
au d\'ebut du chapitre 7 dans \cite{Stil--93}) ~;
\end{itemize}
voir \cite{Dehn--11}, \cite{Dehn--12} et \cite{Dehn--14}.
De plus, ces articles contiennent plusieurs observations
aujourd'hui ``\'evidentes'', mais importantes et \`a l'\'epoque originales.
Citons-en une, de l'introduction de \cite{Dehn--11}~: 
un groupe de pr\'esentation finie peut avoir 
des sous-groupes qui ne sont  pas de type fini.
L'exemple de Dehn, celui d'un groupe libre \`a deux g\'en\'erateurs $s_1, s_2$
et du sous-groupe normal engendr\'e par $s_1$,
est ``\'evident'' si on pense en termes de rev\^etements.
Selon \cite{Magn--78} (page 138),
Dehn fut le premier \`a consigner cette observation.

\medskip

Le m\'elange des genres ci-dessus est r\'ev\'elateur~:
le d\'eveloppement historique de la th\'eorie combinatoire des groupes
est ins\'eparable de celui de la topologie
des vari\'et\'es de basse dimension, 
ce qui veut dire au sens large de dimension au plus $4$,
et tr\`es souvent, strictement, de dimension~$3$.
Plus pr\'ecis\'ement, et comme l'\'ecrit Stallings 
au tout d\'ebut de \cite{Stal--71}~: \par

``The study of three-dimensional manifolds has
often interacted with a certain stream of group theory,
which is concerned with free groups, free products, finite presentations of groups, 
and similar combinatorial matters.
\par
Thus Kneser's fundamental paper 
had latent implications toward Grushko's Theorem.''
\par

Stallings fait allusion
au  th\'eor\`eme de Kneser sur la d\'ecomposition en somme connexe
d'une $3$-vari\'et\'e \cite{Knes--29} 
et au th\'eor\`eme de Grushko sur les d\'ecompositions
en produits libres des groupes de type fini (1940).
Voir aussi ci-dessous l'\'enonc\'e (*), au  chapitre \ref{section7}.

\section{De l'\'etude des n\oe{}uds aux probl\`emes de d\'ecision
\\
Le lemme de Dehn}
\label{section3}

Soit $K$ un \emph{n\oe{}ud}, c'est-\`a-dire une courbe ferm\'ee simple 
diff\'erentiable (ou PL) dans la sph\`ere $\mathbf S^3$.
Notons $V(K)$ un \emph{voisinage tubulaire de $K$,}
c'est-\`a-dire un tore solide $\mathbf S^1 \times \mathbf D^2$
plong\'e dans $\mathbf S^3$ de telle sorte que l'image de 
$\mathbf S^1 \times \{0\}$ co\"{\i}ncide avec $K$,
et $E(K) = \mathbf S^3 \smallsetminus \operatorname{Int}(V(K))$
\emph{l'ext\'erieur de $K$,}
qui est une vari\'et\'e compacte dont le bord 
est le tore $\partial E(K) = \partial V(K)$.
Le \emph{groupe du n\oe{}ud $K$}, not\'e $G_K$,
est le groupe fondamental de $E(K)$.
[On distingue l'ext\'erieur
du \emph{compl\'ement} $\mathbf S^3 \smallsetminus K$,
qui est une vari\'et\'e ouverte sans bord~; 
l'inclusion de l'ext\'erieur dans le compl\'ement est une \'equivalence d'homotopie.]
\par

Un \emph{m\'eridien de $K$} est une courbe ferm\'ee simple $m$
du tore $\partial E(K)$ qui borde un disque dans $V(K)$.
On v\'erifie qu'une telle courbe est unique \`a isotopie pr\`es dans le tore,
et qu'on peut choisir une courbe ferm\'ee simple $q$ dans ce tore
de telle sorte que, apr\`es choix d'orientations, 
les classes de $m$ et $q$ constituent une base de 
$\pi_1(\partial E(K)) \approx H_1((\partial E(K)) \approx\mathbf Z^2$~;
il y a un choix meilleur que d'autres, qui fait l'objet de la d\'efinition suivante.
\par

Notons $G_K^{ab}$ l'ab\'elianis\'e de $G_K$
(ci-dessous not\'e additivement). 
C'est un groupe cyclique infini,
comme on le voit  sur une pr\'esen\-tation 
(de Wirtinger ou de Dehn)  de $G_K$,
ou aussi comme cela r\'esulte\footnote{Le cas particulier
$\pi_1^{ab} \approx H_1$ de l'isomorphisme de Hurewicz
\'etait connu de Poincar\'e \cite{Poin--95}.}
de la dualit\'e de Poincar\'e
\begin{equation*}
G_K^{ab}  \, \approx \,  
H_1(E(K))  \, \approx \, 
H_1(\mathbf S^3 \smallsetminus K)  \, \overset{\text{dualit\'e}}{\approx} \,  
H^1(K) \, \approx \,  
\mathbf Z .
\end{equation*}
Il en r\'esulte d'une part que le groupe des commutateurs de $G_K$
est le noyau du morphisme d'ab\'elianisation
$G_K \longrightarrow G_K^{ab} \approx \mathbf Z$,
et d'autre part que (comme not\'e dans \cite{Dehn--10}, page 115 de \cite{DeSt--87})

\medskip

(3.1) $G_K$ est ab\'elien si et seulement si $G_K \approx \mathbf Z$.

\medskip\noindent
De plus, $G_K^{ab}$ est engendr\'e par la classe $[m]^{ab}$ de $m$.
Soit $q$ comme ci-dessus~;
soit $x \in \mathbf Z$ l'entier tel que $[q]^{ab} = x[m]^{ab}$~;
notons $p$ une courbe ferm\'ee simple
dont la classe co\"{\i}ncide avec $[q]^{ab} - x[m]^{ab} = 0$.
On v\'erifie qu'une telle courbe est unique \`a isotopie pr\`es dans le tore~;
c'est par d\'efinition le \emph{parall\`ele}\footnote{Dans la litt\'erature,
on trouve souvent le mot ``longitude''~;
c'est un non-sens terminologique, 
puisque, chez les g\'eographes,
une longitude est pr\'ecis\'ement la coordonn\'ee d'un m\'eridien~!
(le \emph{m\'eridien de Greenwich} est celui de longitude 0,
la longitude des Diablerets est proche de 7$^o$).
Par ailleurs, notons qu'on peut aussi d\'efinir
$p$ comme l'intersection  $S \cap \partial E(K)$
pour une surface de Seifert $S$ de $K$,
intersection qu'on montre \^etre ind\'ependante \`a isotopie pr\`es
du choix de $S$.
}
de $K$.
Ainsi les courbes $(m,p)$ d\'efinissent (une fois orient\'ees) 
une base de  $\pi_1(\partial E(K)) \approx \mathbf Z^2$,
unique \`a isotopie pr\`es, donc canonique.
[``Unique'' \`a ceci pr\`es qu'on peut toujours
changer les orientations de $m$ et de $p$.
Toutefois, si on souhaite que 
le coefficient d'enlacement de $m$ et $p$ soit $+1$,
un choix d'orientation pour l'un de $m$ et $p$
impose un choix d'orientation pour l'autre.]
Notons encore que,
comme la classe de $p$ dans $H_1(E(K))$ est nulle, 
la classe de $p$ dans $G_K$ est un produit de commutateurs~;
par suite~:

\medskip

(3.2) si $G_K$ est ab\'elien, alors $p = 1$ dans $G_K$.

\medskip\noindent
(La r\'eciproque r\'esulte du lemme de Dehn, voir plus bas.)

\medskip

Un n\oe{}ud $K \subset \mathbf S^3$ est \emph{trivial}
s'il borde un disque plong\'e dans $\mathbf S^3$.
Le groupe d'un n\oe{}ud trivial est cyclique infini engendr\'e par $[m]$,
et le parall\`ele d\'efinit l'\'el\'ement neutre de ce groupe.

\medskip
\noindent
\textbf{Question naturelle.}
\emph{
Soit $K$ un n\oe{}ud de groupe $G_K \approx \mathbf Z$~;
le n\oe{}ud $K$ est-il trivial~?
}

\medskip

Soient $M$ une $3$-vari\'et\'e et $S$ une surface plong\'ee dans $M$.
Un \emph{disque de compression} pour $S$ est
un disque $D$ plong\'e dans $M$ tel que
(i) $D \cap S = \partial D$
et (ii) $\partial D$ ne borde aucun disque dans $S$.
[Nous utilisons la notation $D$ pour \emph{un} disque plong\'e quelque part,
et $\mathbf D^2$ pour \emph{le} disque unit\'e du plan euclidien.]

\medskip
\noindent
\textbf{Lemme facile.}
\emph{
Soient $K$ un n\oe{}ud et $E(K)$, $G_K$ comme plus haut.
S'il existe un disque de compression $D$ 
pour le bord $\partial E(K)$ de $E(K)$,
alors le n\oe{}ud $K$ est trivial.
}

\medskip

\noindent \emph{D\'emonstration.}
On se r\'ef\`ere \`a un m\'eridien $m$ et un parall\`ele $p$ de $K$,
orient\'es, dont les classes constituent une base de 
$H_1(\partial E (K))$.
A priori, le bord de $D$ (orient\'e~!) d\'efinit deux classes d'homologie
\begin{equation*}
\aligned
\left[ \partial D \right]_{\partial E (K)} \, &= \,
x\left[ m \right]_{\partial E (K)} + y\left[ p \right]_{\partial E (K)}
\in H_1(\partial E(K)) ,
\\
\left[ \partial D \right]_{E (K)} \, &= \,
x\left[ m \right]_{E (K)} 
\in H_1(E(K))
\endaligned
\end{equation*}
avec $x,y \in \mathbf Z$.
Or $x=0$ car $\partial D$ est un bord dans $E(K)$,
et $y = \pm 1$ car\footnote{Rappel~:
sur un $2$-tore $T^2$, une courbe ferm\'ee  
de classe d'homologie $x[m] + y[p]  \in H_1(T)$
est isotope \`a une courbe ferm\'ee \emph{simple}
si et seulement si ou bien $x=y=0$ ou bien
$x,y \in \mathbf Z$ sont premiers entre eux.
Voir la section 2.C de \cite{Rolf--76}. 
} 
$\partial D$ est une courbe ferm\'ee simple
dans $\partial E(K)$.
\par

En recollant convenablement \`a $D$ un anneau de bord
$\partial D \cup K$, on obtient un disque plong\'e 
l\'eg\`erement plus grand que $D$, de bord $K$.
Par suite, $K$ est trivial.
\hfill $\square$

\medskip

Revenons \`a la ``question naturelle'' ci-dessus,
pour un n\oe{}ud $K$ de groupe $G_K \approx \mathbf Z$.
On a donc $[p] = 0 \in \mathbf Z$, de sorte qu'il existe
une homotopie $\varphi : \mathbf D^2 \longrightarrow E(K)$
telle que $\varphi(\partial \mathbf D^2) = p$.
Sans perte de g\'en\'eralit\'e, on peut supposer que $\varphi$
est un plongement dans un petit voisinage de $p$,
et que l'\emph{ensemble singulier}
\begin{equation*}
\operatorname{Sing}(\varphi) \, = \, 
\overline{ 
\left\{ x \in \mathbf D^2 \hskip.1cm \vert \hskip.1cm
\vert \varphi^{-1}(\varphi (x)) \vert \ge 2 \right\} 
}
\end{equation*}
est contenu dans l'int\'erieur de $\mathbf D^2$.
Une telle application s'appelle un \emph{disque de Dehn}.
\par

Si l'application $\varphi$ \'etait un plongement,
le n\oe{}ud $K$ serait trivial, vu le ``lemme facile''.
Voil\`a donc de quoi sugg\'erer la formulation
du lemme suivant
(pour la formulation originale, 
voir \cite{Dehn--10}, page 147, 
et \cite{DeSt--87}, page 102).

\medskip
\noindent
\textbf{Lemme de Dehn.}
\emph{
Soient $M$ une $3$-vari\'et\'e et
$\varphi : \mathbf D^2 \longrightarrow M$  
un disque de Dehn.
}
\par
\emph{Alors il existe un plongement ${\psi : \mathbf D^2 \longrightarrow M}$
tel que les restrictions de $\varphi$ et $\psi$  \`a $\partial \mathbf D^2$ co\"incident.
}

\smallskip

\emph{Remarques.}
(i) 
Noter que l'image du bord $\varphi(\partial \mathbf D^2)$ 
est une courbe ferm\'ee \emph{simple} dans  $M$.
Noter aussi que le lemme \emph{ne dit pas} que $\psi$ 
est une petite d\'eformation de $\phi$.
\par
(ii) 
Le cas particulier du lemme qui permet de r\'epondre \`a la
``question naturelle'' formul\'ee plus haut est celui
o\`u la courbe $\varphi(\partial \mathbf D^2)$ est contenue dans le bord 
d'une vari\'et\'e compacte orientable avec bord.
\par
(iii)
Il n'est pas n\'ecessaire de supposer $M$ compacte et orientable.
La r\'eduction du cas non orientable au cas orientable
est due \`a Johansson (1938).
\par
(iv)
Si on sait chercher aussi bien que Gordon, 
on trouve le lemme de Dehn chez Poincar\'e d\'ej\`a, dans \cite{Poin--04}~;
voir la page 475 de \cite{Gord--99}.
\par
(v)
Comme le note Dehn, l'\'enonc\'e a un analogue imm\'ediat
en une dimension de moins~: 
\emph{dans une surface (connexe), 
deux points (qui peuvent \^etre connect\'es par une courbe) 
peuvent \^etre connect\'es par une courbe simple.}

\medskip

Une vingtaine d'ann\'ees plus tard,
plusieurs lecteurs ont d\'ecel\'e une insuffisance s\'erieuse
dans la d\'emonstration du ``lemme de Dehn''~:
Hellmuth Kneser\footnote{
Cette critique de la ``d\'emonstration'' de Dehn appara\^{\i}t \'egalement
dans la  note ajout\'ee aux \'epreuves de \cite{Knes--29}~;
voir aussi \cite{Gord--05}.
Outre ses contributions fondamentales sur les vari\'et\'es de dimension trois,
Hellmuth Kneser (1898-1973)
a des r\'esultats importants dans un tr\`es grand nombre
de domaines math\'ematiques \cite{Knes--05}.
Ne pas le confondre avec son p\`ere Adolf Kneser (1862-1930), 
qui s'illustra dans des travaux d'\'equations diff\'erentielles 
et de calcul des variations,
ni avec son fils Martin Kneser (1928--2004),
connu notamment en combinatoire, formes quadratiques et groupes alg\'ebriques.}, 
dans une lettre \`a Dehn
dat\'ee du 22 avril 1929,
ainsi que van Kampen, Frankl et Pontryagin
(voir \cite{Eppl--99b}, page 279).
Une d\'emonstration fut finalement obtenue par
Papakyriakopoulos, connu sous le diminutif de ``Papa'' \cite{Papa--57b}~;
voir aussi \cite{Stal--01}.

\medskip
\noindent
\textbf{Cons\'equence 
pour les n\oe{}uds.}
\emph{
Soient $K$ un n\oe{}ud, $G_K$ son groupe,
$\iota : \pi_1(\partial E(K)) \longrightarrow G_K$
l'homomorphisme induit par l'inclusion de $\partial E(K)$ dans $E(K)$,
et $p$ un  parall\`ele de $K$,
vu ici comme \'el\'ement de $\pi_1(\partial E(K))$.
Alors les conditions suivantes sont \'equivalentes~:
\begin{itemize}
\item[(i)] 
le n\oe{}ud $K$ est  trivial~;
\item[(ii)] 
le groupe $G_K$ est  ab\'elien~; ou encore, voir (3.1),
\item[]
le groupe $G_K$ est  cyclique infini~;
\item[(iii)]
$\iota(p) = 1$~;
\item[(iv)]
l'homomorphisme $\iota$ n'est pas injectif.
\end{itemize}
}

\medskip\noindent\emph{D\'emonstration.}
Les implications (i) $\Rightarrow$ (ii) et (iii) $\Rightarrow$ (iv)
sont imm\'e\-diates.
L'implication (ii) $\Rightarrow$ (iii) r\'esulte de la d\'efinition de $p$,
selon laquelle  l'image de $p$ dans l'ab\'elianis\'e de $G_K$ est 
l'\'el\'ement neutre.
\par
Les implications non banales sont
(iii) $\Rightarrow$ (i), qui r\'esulte du lemme de Dehn,
et (iv) $\Rightarrow$ (i), qui r\'esulte du th\'eor\`eme du lacet
(voir ci-dessous) et du lemme de Dehn.
\hfill $\square$

\medskip

L'\'equivalence de (i) et (iii) montre bien que,
si on sait d\'ecider entre $\iota(p) = 1 \in G_K$ et $\iota(p) \ne 1$,
on sait d\'ecider de la trivialit\'e ou non de $K$.

\section*{3$^{bis}$. Quelques ``suites'' du lemme de Dehn}

Le ``lemme de Dehn''  \'enonc\'e plus haut fut d\'emontr\'e dans \cite{Papa--57b}.
Ind\'ependamment et juste avant, 
le m\^eme auteur avait aussi montr\'e \cite{Papa--57a}~:

\medskip

\noindent
\textbf{Th\'eor\`eme du lacet.}
\emph{
Si l'homomorphisme $\iota : \pi_1(\partial M) \longrightarrow \pi_1(M)$
induit par l'inclusion n'est pas injectif,
il existe un \'el\'ement non trivial du noyau 
repr\'esent\'e par une courbe ferm\'ee}  simple.

\medskip

En combinant ces deux \'enonc\'es, on obtient le
th\'eor\`eme du disque, parfois nomm\'e ``Loop + Dehn Theorem''~; 
Stallings en fournit une nouvelle d\'emonstration,
ainsi que des g\'en\'eralisations, dans \cite{Stal--60}
(voir aussi dans \cite{Gord--05} un commentaire 
sur cet \'enonc\'e et sur l'article de Kneser \cite{Knes--29}).
Un disque $D$ plong\'e dans $M$ est \emph{proprement plong\'e} 
si $D \cap \partial M = \partial D$.

\medskip

\noindent
\textbf{Cas particulier du th\'eor\`eme du disque de Stallings.}
\emph{Si l'homo\-morphisme $\iota : \pi_1(\partial M) \longrightarrow \pi_1(M)$
induit par l'inclusion n'est pas injectif,
alors il existe un \'el\'ement non trivial du noyau 
repr\'esent\'e par une courbe ferm\'ee simple
bordant un disque proprement plong\'e dans~$M$.
}

\medskip

Dans la foul\'ee, les topologues des ann\'ees 1950 et 1960 ont montr\'e
les r\'esultats suivants (qui ne sont pas \'enonc\'es ici 
de la mani\`ere la plus g\'en\'erale possible). 
En voici le canevas~:
soient $M$ une $3$-vari\'et\'e, $F$ une surface et
$\varphi : F \longrightarrow M$ une application continue 
qui est ``essentielle'', 
c'est-\`a-dire ``non triviale homotopiquement''
(\`a d\'efinir de cas en cas, selon $F$)~;
alors il existe un plongement essentiel de $F$ dans $M$
(assertion \`a prendre parfois \emph{cum grano salis},
voir ci-dessous la ``cons\'equence du th\'eor\`eme du tore'').
Ces r\'esultats montrent que les propri\'et\'es significatives de $\pi_1(M)$ 
refl\`etent des propri\'et\'es g\'eom\'etriques de $M$.

\medskip

\noindent
\textbf{Th\'eor\`eme de la sph\`ere} (\cite{Papa--57b}, J.H.C. Whitehead, 1958).
\emph{
Si une $3$-vari\'et\'e orientable $M$ est telle que $\pi_2(M) \not\approx 0$,
il existe un plongement de $\mathbf S^2$ dans $M$
repr\'esentant un \'el\'ement non nul de $\pi_2(M)$.
}

\medskip

Comme toute $2$-sph\`ere plong\'ee dans un ext\'erieur de n\oe{}ud $E(K)$
borde une $3$-boule (un r\'esultat d'Alexander),
il en r\'esulte que le rev\^ete\-ment universel de $E(K)$ 
est une $3$-vari\'et\'e \emph{contractile}.

\medskip

Notons $\mathbf A$ l'\emph{anneau} $\mathbf S^1 \times [0,1]$
et $\alpha$ un arc proprement plong\'e dans $\mathbf A$ 
connectant les deux  composantes 
$\partial_0 \mathbf A$ et $\partial_1 \mathbf A$ 
de son bord $\partial \mathbf A$.
\par

Une application $\varphi : (\mathbf A, \partial \mathbf A)  
\longrightarrow (M, \partial M)$
est \emph{essentielle} si l'homo\-morphisme induit
$\mathbf Z \approx \pi_1(\mathbf A) \longrightarrow \pi_1(M)$
est injectif
et si $\varphi(\alpha)$ n'est pas homotope relativement \`a ses extr\'emit\'es
\`a un arc dans $\partial M$.

\medskip

\noindent
\textbf{Th\'eor\`eme de l'anneau} (Waldhausen \cite{Wald--69}, \cite{CaFe--76}).
\emph{
Soit $M$ une $3$-vari\'et\'e compacte orientable 
et soit $\varphi : (\mathbf A, \partial \mathbf A)  
\longrightarrow (M, \partial M)$
une application essentielle.}
\par

\emph{Alors il existe un plongement $\psi : \mathbf A \longrightarrow M$ tel que,
pour $j=0$ et $j=1$, les images $\varphi(\partial_j \mathbf A)$
et $\psi(\partial_j \mathbf A)$ soient dans 
la m\^eme composante connexe de $\partial M$.
}
\medskip

Ce th\'eor\`eme permet par exemple d'analyser la situation
o\`u deux lacets dans $\partial M$ d\'efinissent des \'el\'ements conjugu\'es
dans $\pi_1(M)$.

\medskip

Il existe aussi un ``th\'eor\`eme du tore'' pour les applications d'un $2$-tore
dans une vari\'et\'e de dimension $3$, dont nous ne citons que
la cons\'equence suivante (th\'eor\`eme 7 de \cite{Feus--76}).

\medskip

\noindent
\textbf{Cons\'equence du th\'eor\`eme du tore} (\cite{CaFe--76}).
\emph{
Soit $K$ un n\oe{}ud tel que\footnote{Autrement 
dit~: soient $K$ un n\oe{}ud tel qu'il existe  une application ``essentielle''
$\varphi : \mathbf T^2 \longrightarrow E(K)$.} 
$\pi_1(E(K))$ poss\`ede un sous-groupe 
isomorphe \`a $\mathbf Z^2$ qui
n'est pas conjugu\'e \`a un sous-groupe 
du groupe p\'eriph\'erique $\pi_1(\partial E(K))$.
Alors~:
\begin{itemize}
\item[-]
ou bien $K$ est un n\oe{}ud du tore,
\item[-]
ou bien il existe un plongement $\psi$ du $2$-tore $\mathbf T^2$ dans $E(K)$
tel que $\psi_*\left( \pi_1(\mathbf T^2) \right)$ n'est pas conjugu\'e 
\`a un sous-groupe du groupe p\'eriph\'erique $\pi_1(\partial E(K))$,
de sorte que $K$ est un satellite
(par exemple une somme connexe de deux n\oe{}uds non triviaux).
\end{itemize}
}

Rappelons qu'un \emph{n\oe{}ud du tore} est un n\oe{}ud
isotope \`a une courbe ferm\'ee simple contenue dans le tore standard 
de $\mathbf S^3$, et non trivial.
Un tel n\oe{}ud est caract\'eris\'e (\`a l'orientation pr\`es)
par deux entiers $a,b \ge 2$ premiers entre eux~;
son groupe admet alors la pr\'esentation
$\langle s,t \hskip.1cm \vert \hskip.1cm s^a = t^b \rangle$,
et son centre est cyclique infini, engendr\'e par $s^a$.  
Pour une description du sous-groupe p\'eriph\'erique,
voir par exemple la proposition 3.28 de \cite{BuZi--85}.

\section*{3$^{ter}$. Quelques exemples de ``probl\`emes (ir)r\'esolubles''}

D\'ecrivons au moins un ``processus de d\'ecision''
montrant qu'il existe des groupes \`a probl\`eme du mot r\'esoluble.
Consid\'erons un entier $n \ge 1$ et le \emph{groupe libre} $F_n$ de rang $n$,
avec sa pr\'esentation $\langle s_1, \hdots, s_n \hskip.1cm \vert \hskip.1cm \rangle$.
(Il serait bien s\^ur possible, mais un peu plus compliqu\'e,
de consid\'erer une pr\'esentation finie \emph{arbitraire} de $F_n$.)

Ce processus s'applique \`a un mot $w$ de longueur $\vert w \vert$ 
en les $s_i$ et leurs inverses, et fonctionne en trois pas.
\begin{itemize}
\item[(i)]
Si $\vert w \vert \le 1$, aller en (iii).
Si $\vert w \vert \ge 2$, aller en (ii).
\item[(ii)]
S'il existe dans $w$ deux lettres cons\'ecutives $s_is_i^{-1}$ ou  $s_i^{-1}s_i$,
les supprimer et retourner en (i) avec le mot ainsi obtenu, 
de longueur  $\vert w \vert - 2$. Sinon, aller en (iii).
\item[(iii)]
Si le mot obtenu est vide, \'ecrire $w =_G  1$ et s'arr\^eter.
Sinon, \'ecrire $w \neq_G 1$ et s'arr\^eter.
\end{itemize}
(Il y a une \'etroite analogie avec
l'\emph{algorithme de Dehn} pour les groupes hyperboliques,
d\'ecrit  au chapitre \ref{section6}.)

Pour la description d'un exemple moins imm\'ediat 
mais n\'eanmoins accessible au non-sp\'ecialiste, 
qui est une solution du probl\`eme du mot pour
un \emph{groupe r\'esiduellement fini de pr\'esentation finie,}
nous renvoyons le lecteur au th\'eor\`eme 2.2.5 de \cite{Robi--96}.

\medskip

La d\'emonstration de l'existence de groupes de pr\'esentation finie
\`a probl\`eme du mot non r\'esoluble comprend deux \'etapes.
\begin{itemize}
\item[(SG)]
\emph{Il existe un semi-groupe de pr\'esentation finie \`a probl\`eme du mot
non r\'esoluble.} 
C'est un r\'esultat obtenu ind\'ependamment par
A.A.~Markov\footnote{Il s'agit de Andrei Andreyevich Markov (1903--1979),
fils de  Andrei Andreyevich Markov (1856-1922),
dont le directeur de th\`ese (1884) \'etait Pafnuty Lvovich Chebyshev.
Les ``cha\^{\i}nes de Markov'' portent le nom du p\`ere.} 
et Emil Post (1947)~; sa d\'emonstration utilise
l'existence dans $\mathbf N$ d'un sous-ensemble
``r\'ecursivement \'enum\'e\-rable mais non \'enum\'erable'',
existence qu'on montre en uti\-lisant entre autres
la th\'eorie des \emph{machines de Turing} 
et une variante de l'argument diagonal de Cantor
(celui qui montre que $\mathbf R$ n'est pas d\'enombrable).
\item[(G)]
\emph{R\'eduction du cas des groupes au cas des semi-groupes.}
Il s'agit de montrer qu'il existe un groupe $G$ tel que,
si son probl\`eme du mot \'etait r\'esoluble,
alors le probl\`eme du mot serait aussi r\'esoluble
pour un semi-groupe \`a la Markov-Post,
ce qui n'est pas vrai par (SG).
Cette \'etape, due (autant que je sache  ind\'ependam\-ment)
\`a P.S.~Novikov, W.W.~Boone et J.L.~Britton,
date de la fin des ann\'ees 1950.
Elle  a \'et\'e notoirement simplifi\'ee depuis,
notamment gr\^ace \`a un usage important de la notion
d'extension HNN 
(notion datant d'un article de 1949 par
 Graham Higman, Bernhard Hermann Neumann et Hanna Neumann).
\end{itemize}

Citons un r\'esultat marquant, qui est une caract\'erisation alg\'ebrique
de la solubilit\'e du probl\`eme des mots.

\medskip

\noindent
\textbf{Th\'eor\`eme de Boone-Higman} (1974).
\emph{Un groupe de type fini a un probl\`eme du mot r\'esoluble
si et seulement s'il est isomorphe \`a un sous-groupe
d'un sous-groupe simple d'un groupe de pr\'esentation finie.}

\medskip

Pour une exposition de ce qui pr\'ec\`ede, voir le chapitre 12 de \cite{Rotm--95}.

\medskip
\begin{center}
* * * * * * * * * * * * *
\end{center}
\smallskip

Au d\'ebut des ann\'ees 1930, Magnus a montr\'e que
le probl\`eme du mot est r\'esoluble pour les \emph{groupes \`a un relateur}~;
c'est une cons\'equence du ``Freiheitssatz''
(voir le chapitre II.5 de \cite{ChMa--82}).

Le probl\`eme du  mot pour les \emph{groupes de n\oe{}ud} est r\'esoluble~:
c'est un (cas particulier d'un) r\'esultat de  Waldhausen \cite{Wald--68}~;
en fait, depuis les travaux de Perelman \'etablissant
la conjecture de g\'eom\'etrisation de Thurston
(voir par exemple \cite{BBBMP}), 
on sait que cela vaut pour tous les \emph{groupes fondamentaux 
de  $3$-vari\'et\'es} (voir le chapitre 12 de \cite{Ep+5--92} et
la discussion autour du th\'eor\`eme 3.3.1 de \cite{Brid--02}).
Mieux, la croissance de la fonction de Dehn d'un tel groupe
est au plus exponentielle \cite{Brid--93}.
Le probl\`eme de conjugaison pour les groupes fondamentaux
de $3$-vari\'et\'es est r\'esoluble~: voir \cite{Prea--06} et \cite{Prea}.

Le probl\`eme du  mot pour les \emph{groupes lin\'eaires de type fini} est r\'esoluble.
Le probl\`eme de conjugaison pour les groupes lin\'eaires de pr\'esentation finie 
est r\'esoluble.
Voir les th\'eor\`emes 5.1 et 5.3 de \cite{Mill--92}.
Pour un \emph{groupe de type fini r\'esiduellement fini}, 
le probl\`eme du mot est r\'esoluble dans le cas de pr\'esentation finie
(comme nous y avons d\'ej\`a fait allusion un peu plus haut),
mais pas en g\'en\'eral \cite{Mesk--74}.
Notons que, depuis les succ\`es de Perelman,
on sait que tout groupe fondamental de $3$-vari\'et\'e compacte
est r\'esiduellement fini \cite{Hemp--87}~;
mais on ne sait toujours pas si un tel groupe est lin\'eaire 
(voir n\'eanmoins \cite{AsFr}).
\par

Pour les \emph{groupes hyperboliques au sens de Gromov},
le probl\`eme du mot et le probl\`eme de conjugaison 
sont r\'esolubles \cite{Grom--87}~;
voir aussi \cite{BrHa--99}, pages 448 et suivantes,
ainsi que la fin du chapitre \ref{section6} ci-dessous.
Le probl\`eme d'isomorphisme pour les groupes hyperboliques 
est aussi r\'esoluble~;
voir \cite{Sela--95}, avec hypoth\`eses suppl\'ementaires,
et \cite{DaGu}, pour le cas g\'en\'eral.

Il y a de bonnes raisons d'\'etudier le probl\`eme du mot pour certains groupes
qui ne sont pas de type fini~; noter qu'un mot en un syst\`eme (infini) $S$
de g\'en\'erateurs d'un groupe $G$ est toujours dans un sous-groupe de type fini de $G$.
Par exemple, Lyndon a montr\'e que le probl\`eme du mot est r\'esoluble
dans ``son''  groupe $F^{\mathbf Z [x]}$, 
qui est un groupe tel que ses sous-groupes de type fini 
sont exactement  les ``groupes limites'' de Sela
\cite{Lynd--60}.

\section{Diagrammes et graphes de Cayley, 
\\ 
Dehn's Gruppenbilder}
\label{section4}

Soit $G = \langle S \rangle$ un groupe donn\'e avec
un ensemble fini de g\'en\'erateurs
$S = \{ s_1, s_2, \hdots, s_n \}$.
Le \emph{diagramme de Cayley} ou ``\emph{Dehn Gruppenbild}'' associ\'e
est le graphe orient\'e $\mathcal C (G,S)$
ayant $G$ pour ensemble de sommets,
et dans lequel il y a une ar\^ete orient\'ee de $g$ \`a $h$
et \'etiquet\'ee  $s$ lorsque $g^{-1}h = s \in S \cup S^{-1}$. 
Relevons les points suivants.
\begin{itemize}
\item[-]
Dans $\mathcal C (G,S)$, avec toute ar\^ete de $g$ vers $h$
\'etiquet\'ee $s$, on trouve
une ar\^ete de $h$ vers $g$ \'etiquet\'ee $s^{-1}$.
Toutefois, dans les dessins, on ne repr\'esente souvent qu'une
des deux ar\^etes d'une m\^eme paire, 
celle dont l'\'etiquette $s$ est dans $S$~;
de plus, dans le cas sp\'ecial o\`u $s^2 = 1$ dans $G$
(ce qu'on \'ecrit $s^2 =_G 1$),
on omet souvent d'indiquer l'orientation.
\item[-]
Il peut exister plusieurs ar\^etes de $g$ vers $h$, d'\'etiquettes diff\'erentes.
En effet, il faut comprendre que $S$ n'est pas dans $G$,
mais est une base d'un groupe libre $F$ donn\'e avec
une surjection $\pi : F \longrightarrow G$,
et la restriction de $\pi$ \`a $S$ n'est pas n\'ecessairement injective.
De m\^eme, s'il existe $s \in S$ tel que $\pi(s) = 1$,
il y a des boucles aux sommets de $\mathcal C (G,S)$.
\item[-]
Exemple d'une pr\'esentation dont le diagramme de Cayley 
contient des ar\^etes multiples et des boucles~:
\begin{equation*}
\mathbf Z \, = \, 
\langle s,t,u \hskip.1cm \vert \hskip.1cm  s = t, \hskip.1cm u = 1 \rangle .
\end{equation*}
\item[-]
Le groupe $G$ op\`ere sur $\mathcal C (G,S)$ par automorphismes
de graphe dirig\'e \'etiquet\'e~; 
$x \in G$ applique l'ar\^ete de $g$ \`a $gs$
sur l'ar\^ete de $xg$ \`a $xgs$.
Cette action est simplement transitive.
\end{itemize}
Dans ce texte, nous conviendrons d'appeler \emph{graphe de Cayley}
le graphe g\'eom\'etrique associ\'e au diagramme de Cayley,
c'est-\`a-dire le graphe obtenu en oubliant
les orientations et les \'etiquettes, et en rempla\c cant
chaque couple  $\big((g,gs), (gs, (gs)s^{-1})\big)$ d'ar\^etes orient\'ees
par une seule ar\^ete $\{g,gs\}$.
\par
Autres points de vue~:
le graphe de Cayley est un espace m\'etrique d\'enombrable
qui est connexe par  pas de longueur $1$,
c'est aussi le $1$-squelette d'un $2$-complexe convenable
(voir plus bas au chapitre  \ref{section7}).

\subsection{Un  diagramme de Cayley
d'un groupe d'isom\'etries propres
du plan hyperbolique}
\label{ 1}

Le premier exemple donnant lieu \`a une figure dans \cite{Dehn--10}
est celui associ\'e \`a la pr\'esentation
\begin{equation*}
\langle s_1, s_2 \hskip.1cm \vert \hskip.1cm (s_1)^3, (s_2)^5, (s_1s_2)^2 \rangle
\end{equation*}
du groupe altern\'e $A_5$ d'ordre $60$.
Les notes \cite{Dehn--a} contiennent d'autres exemples, 
dans l'ordre~:
les groupes sym\'etriques $\Sigma_3$ et $\Sigma_4$, 
les groupes altern\'es $A_4, A_5$,
le groupe $A_5 \times \{1,-1\}$ des isom\'etries laissant invariant
un icosa\`edre r\'egulier,
et le groupe sym\'etrique
\begin{equation*}
\Sigma_5 \, = \, \left\langle s_1, s_2, s_3 \hskip.2cm \Bigg\vert \hskip.2cm 
\aligned
&s_1^2 = s_2^5 = (s_1s_2)^2 = 1
\\
&s_3^2 = 1, \hskip.2cm 
s_3^{-1} s_1 s_3 = s_1, \hskip.2cm
s_3^{-1} s_2 s_3 = s_2
\endaligned
\right\rangle .
\end{equation*}
Dehn remarque que 
$s_2 = (1 \hskip.1cm 2 \hskip.1cm 3 \hskip.1cm 4 \hskip.1cm 3)$
et $s_3 = (1 \hskip.1cm 2)$
suffisent \`a engendrer $\Sigma_5$, 
et que celui-ci est un quotient du groupe
\begin{equation*}
G \, = \, \langle s_2, s_3 \hskip.1cm \vert \hskip.1cm
(s_2)^5, (s_3)^2, (s_3s_2)^4 \rangle .
\tag{$P_{5,4}$}
\end{equation*}
Par un argument g\'eom\'etrique, 
il montre que $G$ est un groupe infini, 
et il en d\'ecrit le graphe de Cayley relativement \`a $\{s_2,s_3\}$~;
nous comprenons l'argument de Dehn en le reformulant comme suit.

Consid\'erons dans $\mathcal H^2$ un triangle isoc\`ele $T$
ayant un sommet $A$ d'angle $2\pi/5$ 
et deux sommets $B,C$ d'angle $2\pi/8$.
Soient $s_2$ une rotation de centre $A$ et d'angle $2\pi/5$,
et $s_3$ un demi-tour centr\'e au milieu du segment $BC$.
Il r\'esulte d'un th\'eor\`eme de Poincar\'e 
sur les polygones g\'en\'erateurs des groupes fuchsiens
que le groupe d'isom\'etries de $\mathcal H^2$ engendr\'e par $s_2$ et $s_3$ 
d'une part a $T$ pour domaine fondamental
et d'autre part a la pr\'esentation ($P_{5,4}$).
(Pour le th\'eor\`eme de Poincar\'e, que Dehn n'invoque pas explicitement,
voir \cite{Poin--82}, ou l'exposition de \cite{Iver--92}.)
Le graphe de Cayley correspondant est le graphe dual
du pavage $\left( gT \right)_{g \in G}$ de $\mathcal H^2$~; autrement dit~:
\begin{itemize}
\item[]\emph{le graphe de Cayley associ\'e \`a la pr\'esentation ($P_{5,4}$)
est un graphe trivalent, qui
est  le $1$-squelette d'un pavage de $\mathcal H^2$ par des
pentagones r\'eguliers (centr\'es aux points de la forme $gA$)
et des octogones r\'eguliers (centr\'es aux points de la forme $gB$ et $gC$),
avec un pentagone et deux octogones incidents \`a chaque sommet.
}
\end{itemize}

Toujours dans \cite{Dehn--a},  
Dehn ajoute que $G$ s'ins\`ere dans la famille
\begin{equation*}
G_{\alpha,\beta} \, = \, \langle s_1, s_2 \hskip.1cm \vert \hskip.1cm
(s_1)^\alpha, (s_2)^2, (s_2s_1)^\beta \rangle ,
\tag{$P_{\alpha,\beta}$}
\end{equation*}
o\`u $\alpha$ et $\beta$ sont deux entiers.
Le groupe $G_{\alpha,\beta}$ est infini si 
\begin{equation*}
\frac{\alpha - 2}{\alpha} \hskip.1cm \pi 
+ 2 \hskip.1cm \frac{2 \beta - 2}{2 \beta} \hskip.1cm \pi
\, \ge \, 2 \pi ,
\end{equation*}
et les tuiles du pavage correspondant de $\mathcal H^2$
sont alors des $\alpha$-gones r\'eguliers
et des $2 \beta$-gones r\'eguliers.
\par

Ajoutons encore une remarque.
Soient $s'_3$ une sym\'etrie de $\mathcal H^2$ d'axe prolongeant le c\^ot\'e $BC$,
et $G'$ le groupe d'isom\'etries de $\mathcal H^2$ engendr\'e par $s_2$ et $s'_3$.
Le th\'eor\`eme de Poincar\'e fournit la pr\'esentation
\begin{equation*}
G' \, = \, \langle s_2, s'_3 \hskip.1cm \vert \hskip.1cm
(s_2)^5, (s'_3)^2, (s'_3{s_2}^{-1}s'_3s_2)^2 \rangle .
\tag{$P_{5,4}'$}
\end{equation*}
Comme $T$ est \'egalement domaine fondamental 
pour l'action de $G'$ sur $\mathcal H^2$,
le graphe de Cayley de ($P_{5,4}'$)
est isom\'etrique \`a celui de ($P_{5,4}$).
Notons que le groupe $G'$, dont la pr\'esentation ($P_{5,4}'$)
montre imm\'ediatement que son ab\'elianis\'e est cyclique d'ordre $10$,
n'est pas isomorphe au groupe $G$, dont l'ab\'elianis\'e est d'ordre $2$.

\subsection{G\'en\'eralit\'es sur les diagrammes de Cayley}
\label{subsection4.2}

Il est naturel de chercher des conditions n\'ecessaires et/ou suffisantes
pour qu'un diagramme (graphe orient\'e et \'etiquet\'e)
soit le diagramme de Cayley associ\'e 
\`a une pr\'esentation de groupe.
Avant de formuler un \'enonc\'e, rappelons ceci.

Consid\'erons un \emph{diagramme r\'egulier}\footnote{Attention~:
ce n'est pas l'usage du mot ``r\'egulier'' \`a la page 62 de \cite{MaKS--66},
qui me semble propice aux confusions avec 
la signification du mot pour les graphes.
Nous pr\'ef\'erons suivre \cite{Cann--02} et
traduire le ``regular'' de Magnus-Karras-Solitar par ``homog\`ene''.},
c'est \`a dire un graphe orient\'e, 
\'etiquet\'e par un ensemble de la forme $\{s_1, s_1^{-1}, \hdots, s_n,s_n^{-1}\}$,
chaque sommet \'etant incident \`a 
exactement une ar\^ete entrante \'etiquet\'ee $s_i^{-1}$
et une ar\^ete sortante \'etiquet\'ee $s_i$, pour $i = 1,\hdots,n$,
et chaque ar\^ete \'etiquet\'ee $s_i$ d'un sommet $x$ vers un sommet $y$
\'etant appari\'ee \`a une ar\^ete \'etiquet\'ee $s_i^{-1}$
du sommet $y$ vers le sommet $x$.
(A ceci pr\`es que, si $s_i^2 =_G 1$,
chaque ar\^ete  \'etiquet\'ee $s_i$ de $x$ vers $y$
est appari\'ee \`a une ar\^ete \'etiquet\'ee $s_i$ de $y$ vers $x$.)
Etant donn\'e un tel diagramme,
on associe naturellement \`a tout chemin dans ce graphe
un mot en les $s_i$ et leurs inverses.
Le diagramme est alors \emph{homog\`ene} si,
pour tout couple de chemins d'origines diff\'erentes
et de m\^eme mot associ\'e, les chemins sont en m\^eme temps
ferm\'es ou non.

\medskip

\emph{Exemples de diagrammes r\'eguliers non homog\`enes.}
(i) Soit $D$ un diagramme \`a trois sommets $x,y,z$
avec une boucle \'etiquet\'ee $s_1$ en $x$,
une ar\^ete \'etiquet\'ee $s_2$ de $x$ vers $y$ et une autre de $y$ vers $x$,
une ar\^ete \'etiquet\'ee $s_1$ de $y$ vers $z$ et une autre de $z$ vers $y$,
et une boucle \'etiquet\'ee $s_2$ en $z$
(figure 4, page 62 de \cite{MaKS--66}).
Alors $D$ est r\'egulier, mais non homog\`ene
car le chemin d'\'etiquette $s_1$ issu de $x$ est ferm\'e
et celui de m\^eme \'etiquette issu de $y$ ne l'est pas.
\par
(ii) Soit $P$ le graphe de Petersen, qu'on peut d\'efinir
comme le graphe dont les sommets sont les 
sous-ensembles \`a deux \'el\'ements de l'ensemble $\{1,2,3,4,5\}$,
avec une ar\^ete entre $\{i,j\}$ et $\{k,\ell\}$
lorsque $i,j,k,\ell$ sont tous distincts.
C'est un graphe r\'egulier de degr\'e trois,
\`a dix sommets, de groupe de sym\'etrie isomorphe
au groupe sym\'etrique  $\Sigma_5$.
Si $P$ \'etait le graphe sous-jacent \`a un diagramme de Cayley
$\mathcal C (G,S)$, l'ensemble g\'en\'erateur $S$
contiendrait au moins un \'el\'ement d'ordre $2$
(pour obtenir un graphe de degr\'e impair)~;
mais ceci est impossible car tout \'el\'ement d'ordre $2$
du groupe de sym\'etrie $\Sigma_5$ de $P$ poss\`ede des sommets fixes.

\medskip

\noindent
\textbf{Proposition.}
\emph{Soit $\mathcal C$ un diagramme connexe r\'egulier, 
\'etiquet\'e par un ensemble de la forme 
$S = \{s_1, s_1^{-1}, \hdots, s_n,s_n^{-1}\}$, 
et soit $G$ son groupe d'automorphismes
(en tant que graphe dirig\'e \'etiquet\'e).}
\par
\emph{Alors $\mathcal C$ est un diagramme de Cayley 
si et seulement s'il est homog\`ene.}
\par
\emph{Lorsque c'est le cas, c'est le diagramme de Cayley du groupe $G$
rela\-tivement au syst\`eme de g\'en\'erateurs $S$.
En particulier, l'action de $G$ sur l'ensemble des sommets de $\mathcal C$
est simplement transitive.}

\medskip

\noindent \emph{R\'ef\'erences pour la d\'emonstration.}
Voir le th\'eor\`eme 1.6, page 63,
ainsi que l'exercice 15, page 69, de \cite{MaKS--66}.
\par
On trouve aussi une bonne discussion sur les graphes de Cayley
dans l'article de Cannon \cite{Cann--02}.
\hfill $\square$

\medskip

R\'eciproquement, soit $G$ un groupe donn\'e 
par une pr\'esentation finie $\langle S \hskip.1cm \vert \hskip.1cm R \rangle$
et soit $\mathcal C$ un diagramme connexe r\'egulier,
\'etiquet\'e par $S \cup S^{-1}$, homog\`ene,
et dans lequel tout chemin \'etiquet\'e par une relation de $R$ est ferm\'e.
Alors $\mathcal C$ est le diagramme de Cayley \emph{d'un quotient} de $G$ mais,
en g\'en\'eral, pas de $G$ lui-m\^eme.
L'exemple le plus simple est sans doute celui d'un circuit de longueur $n$
convenablement orient\'e et \'etiquet\'e,
qui est le diagramme de Cayley d'un quotient fini
$\mathbf Z / n\mathbf Z = \langle s \hskip.1cm \vert \hskip.1cm s^n \rangle$
du groupe cyclique infini
 $\mathbf Z = \langle s \hskip.1cm \vert \hskip.1cm \rangle$.

Pour pouvoir d\'ecider si $\mathcal C$ est vraiment 
le diagramme $\mathcal C  (G,S)$,
il faut par exemple conna\^{\i}tre une forme normale des \'el\'ements de $G$
\'ecrits en termes des g\'en\'erateurs de $S$,
ou encore avoir r\'esolu le probl\`eme du mot pour $G$.
Autrement dit~:
\begin{itemize}
\item[] \emph{dessiner des parties arbitrairement grandes}
\item[] \emph{du  ``Dehn Gruppenbild'' de 
$G = \langle S \hskip.1cm \vert \hskip.1cm R \rangle$}
\item[] \emph{\'equivaut \`a r\'esoudre le probl\`eme du mot pour $G$.}
\end{itemize}
Davantage sur cette \'equivalence dans \cite{Cann--02}.

\subsection{Quelques diagrammes pour les groupes cycliques
et certains de leurs produits libres}
\label{subsection4.3}

Ce num\'ero est un ``exercice d'\'echauf\-fement''
en vue du num\'ero \ref{subsection4.5}.
\par
Consid\'erons le groupe $\mathbf Z$ des entiers rationnels,
engendr\'e par $\{1,k\}$, avec $k \ge 2$,
et le diagramme de Cayley $\mathcal C (\mathbf Z, \{1,k\})$ correspondant.
\par

Si $k=2$, ce diagramme est une \'echelle infinie,
avec les entiers pairs sur un montant infini 
aux altitudes $\hdots, -4, -2, 0, 2, 4, \hdots$,
les entiers impairs sur l'autre montant infini
aux altitudes $\hdots, -3, -1, 1, 3, \hdots$,
tous les segments de ces montants dirig\'es vers le haut
et \'etiquet\'es par $2$,
et des \'echelons obliques liant $n$ \`a $n+1$
\'etiquet\'es par $1$.
Le dessin obtenu est la figure 11 de \cite{Dehn--10}.
\par

Si $k=3$, le diagramme peut \^etre d\'ecrit comme
un prisme triangulaire,
avec trois montants infinis contenant les sommets
qui correspondent aux entiers congrus \`a $0,1,2$
respectivement, et une spirale montante
lin\'eaire par morceaux, chaque morceau
connectant un entier $n$ sur l'un des montants
\`a l'entier $n+1$ sur le prochain montant.
\par

De m\^eme, pour tout $k \ge 3$, 
le diagramme de Cayley $\mathcal C (\mathbf Z, \{1,k\})$
peut \^etre vu comme un prisme
infini sur un $k$-gone r\'egulier
avec une spirale lin\'eaire par morceaux.

\medskip

Consid\'erons ensuite deux entiers $k,\ell \ge 2$,
un groupe cyclique $C_k$ d'ordre $k$ engendr\'e par un g\'en\'erateur $\mu$,
un groupe cyclique $C_\ell$ d'ordre $\ell$ engendr\'e par un g\'en\'erateur $\nu$,
et leur produit libre $C_k \ast C_\ell$ engendr\'e par $\{\mu, \nu\}$.
\par

Le diagramme de Cayley $\mathcal C (C_k \ast C_\ell , \{\mu,\nu\})$
peut \^etre dessin\'e comme suit, par couches~:
la premi\`ere couche est un $k$-gone r\'egulier,
dont les sommets correspondent aux \'el\'ements 
$\mu^m, m = 0,1,\hdots,k-1$, du groupe $C_k \ast C_\ell$~;
la seconde couche est une couronne de $\ell$-gones r\'eguliers,
au nombre de $k$, 
chacun ayant exactement un sommet commun 
avec le $k$-gone de la couche pr\'ec\'edente~;
la couche suivante est une couronne de $k$-gones r\'eguliers,
au nombre de $k(\ell - 1)$, 
chacun ayant exatement un sommet commun 
avec un $\ell$-gone de la couche pr\'ec\'edente~; etc.
(On convient qu'un $2$-gone r\'egulier est un segment.)

``Vu d'assez loin'', $\mathcal C (C_k \ast C_\ell , \{\mu,\nu\})$
se rapproche d'un arbre biparti ayant alternativement
des sommets de degr\'e $k$ et des sommets de degr\'e $\ell$.
En termes plus techniques, $\mathcal C (C_k \ast C_\ell , \{\mu,\nu\})$
est quasi-isom\'etrique \`a un tel arbre.

\subsection{Le  diagramme de \cite{Dehn--10}
pour le groupe du n\oe{}ud de tr\`efle}
\label{subsection4.4}
Soient $K$ un n\oe{}ud de tr\`efle et $G$ son groupe.
L'un des diagrammes usuels  repr\'esentant $K$
fournit la pr\'esentation de Dehn\begin{equation*}
G \, = \, 
\langle a,b,c,d \hskip.1cm \vert \hskip.1cm
ad^{-1}b = bd^{-1} c = cd^{-1}a = 1 \rangle ,
\tag{$P_1$}
\end{equation*}
apparaissant au \S~5 du chapitre II de \cite{Dehn--10}.
La pr\'esentation $(P_1)$ de $G$ ayant $4$ g\'en\'erateurs
dont aucun n'est d'ordre $2$ dans $G$,
le graphe de Cayley  $\mathcal C (G, \{a,b,c,d\})$
est  r\'egulier de degr\'e $8$~;
d\'ecrivons-le. 
\par

Soit d'abord $L$ l'\'echelle obtenue \`a partir du graphe
$\mathcal C (\mathbf Z, \{1,2\})$ du num\'ero  pr\'ec\'edent
en \'etiquetant toutes les ar\^etes verticales des deux montants infinis par $d$
et les ar\^etes obliques montantes du zig-zag m\'edian 
par une suite p\'eriodique $\hdots, c, b, a, c, b, a, \hdots$
(la lettre $L$ se r\'ef\`ere \`a ``Leiter = \'echelle'').
Soit par ailleurs $B$ l'arbre trivalent r\'egulier
(la lettre $B$ se r\'ef\`ere \`a ``Baum = arbre'').
Alors
\begin{itemize}
\item[]\emph{le diagramme $\mathcal C (G, \{a,b,c,d\})$ 
de la pr\'esentation ($P_1$)
est la r\'eunion d'une famille $\left( L_e \right)_{e \in  E (B)}$ 
de copies de $L$ index\'ee par l'ensemble $E (B)$ 
des ar\^etes g\'eom\'etriques de $B$,
chaque copie de $L$ \'etant recoll\'ee \`a deux autres copies
sur chacun des montants verticaux (donc \`a quatre autres copies en tout),
de telle sorte que chaque sommet du graphe r\'esultant
soit incident \`a quatre ar\^etes entrantes, d'\'etiquettes $a,b,c,d$,
et quatre ar\^etes sortantes, d'\'etiquettes $a,b,c,d$.}
\end{itemize}
Comme d\'ej\`a dit, il convient de ne pas dessiner
les ar\^etes \'etiquet\'ees par l'une des lettres
$a^{-1}, b^{-1}, c^{-1}, d^{-1}$.

\medskip

Le diagramme ainsi d\'ecrit est a priori
le diagramme de Cayley d'un quotient (sic !) de $G$.
Pour s'assurer que c'est vraiment le diagramme de Cayley de $G$,
il faut se convaincre que tout chemin ferm\'e dans ce graphe
peut \^etre homotop\'e \`a un chemin constant
via des ``homotopies \'el\'ementaires'', chacune correspondant
\`a un triangle associ\'e \`a l'une des relations de la pr\'esentation 
($P_1$)~; le fait que le diagramme ``se projette''
sur l'arbre $B$ joue un r\^ole important dans l'argument .....
Voir \cite{Dehn--10}, page 117 de \cite{DeSt--87},
compl\'et\'e par \cite{Dehn--11}, pages 173-174 de \cite{DeSt--87}.

\subsection{Des diagrammes de Cayley
pour les groupes des n\oe{}uds toriques}
\label{subsection4.5}

Nous allons d\'ecrire une autre pr\'esentation de $G$
et son diagramme de Cayley.
La pr\'esentation ($P_1$) peut  s'\'ecrire
\begin{equation*}
G \, = \, 
\langle a,b,c,d \hskip.1cm \vert \hskip.1cm
ba =cb = ac = d \rangle ,
\tag{$P'_1$}
\end{equation*}
de sorte qu'on a aussi
\begin{equation*}
G \, = \, 
\langle a,b,c \hskip.1cm \vert \hskip.1cm
ba =cb = ac  \rangle ,
\tag{$P_2$}
\end{equation*}
qui est d'ailleurs  une pr\'esentation de Wirtinger de $G$,
et encore
\begin{equation*}
G \, = \, 
\langle a,b \hskip.1cm \vert \hskip.1cm
bab =aba  \rangle ,
\tag{$P_3$}
\end{equation*}
montrant que $G$ est isomorphe au
``groupe d'Artin des tresses \`a trois brins''.
En posant  $v = bab$, on trouve
$G  = \langle d,v \hskip.1cm \vert \hskip.1cm
d^3 = v^2  \rangle$,
ou encore
\begin{equation*}
G \, = \, 
\langle t,u,v \hskip.1cm \vert \hskip.1cm
t = u^3 = v^2  \rangle .
\tag{$P_4$}
\end{equation*}
Cette derni\`ere pr\'esentation montre
que $G$ s'ins\`ere dans une suite exacte de la forme
\begin{equation*}
1 
\hskip.2cm \longrightarrow \hskip.2cm
\mathbf Z = \langle t \rangle  
\hskip.2cm \longrightarrow \hskip.2cm
G = \langle t,u,v \rangle  
\hskip.2cm  \overset\pi\longrightarrow \hskip.2cm
C_3 \ast C_2 
\hskip.2cm \longrightarrow \hskip.2cm 
1 ,
\end{equation*}
o\`u $\mathbf Z = \langle t \rangle$ d\'esigne le centre de $G$,
engendr\'e par $t$,
et o\`u $C_3 \ast C_2$ d\'esigne comme plus haut le produit libre
d'un groupe cyclique d'ordre~$3$, 
engendr\'e par la classe  de $u$ modulo le centre de $G$,
et d'un groupe cyclique d'ordre~$2$, engendr\'e par celle de $v$.

Les \'el\'ements de 
$C_3 \ast C_2 = \langle \mu,  \nu \hskip.1cm \vert \hskip.1cm
\mu^3 = \nu^2 = 1 \rangle$
sont en bijection avec les mots de la forme
\begin{equation*}
\mu^{\epsilon_0}
\nu \mu^{\epsilon_1} 
\nu \mu^{\epsilon_2} 
\cdots
\nu \mu^{\epsilon_n}
\end{equation*}
o\`u $n \ge 0$, $\epsilon_0, \epsilon_n \in \{0,1,2\}$,
$\epsilon_1, \hdots, \epsilon_{n-1} \in \{1,2\}$.
Il en r\'esulte que les \'el\'ements de $G$
sont en bijection avec les mots de la forme
\begin{equation*}
t^k
u^{\epsilon_0} v u^{\epsilon_1} v u^{\epsilon_2} \cdots
v u^{\epsilon_n}
\tag{N}
\end{equation*}
o\`u $k \in \mathbf Z$ et $n, \epsilon_0, \hdots, \epsilon_n$
comme ci-dessus.
Le diagramme de Cayley $\mathcal C (G,\{t,u,v\})$, de degr\'e $6$,
se projette donc naturellement sur le diagramme
$\mathcal C ( C_3 \ast C_2, \{\mu,\nu\})$, de degr\'e $3$,
d\'ecrit selon ses couches au num\'ero pr\'ec\'edent.
\par

On peut d'ailleurs g\'en\'eraliser sans peine
au cas de deux entiers $k,\ell \ge 2$
et du groupe
\begin{equation*}
G_{k,\ell} \, = \, 
\langle t,u,v \hskip.1cm \vert \hskip.1cm
t = u^k = v^\ell  \rangle ;
\tag{$P_{k,\ell}$}
\end{equation*}
lorsque de plus $k$ et $\ell$ sont premiers entre eux,
$G_{k,\ell}$ est le groupe d'un \emph{n\oe{}ud du tore}.
Le diagramme $\mathcal C (G_{k,\ell},\{t,u,v\})$ peut \^etre d\'ecrit comme suit~:
au-dessus de la premi\`ere couche du diagramme
de $\mathcal C (C_k \ast C_\ell, \{\mu,\nu\})$
d\'ecrit en seconde partie du num\'ero 
\ref{subsection4.3}, 
un prisme infini sur un $k$-gone r\'egulier
comme en premi\`ere partie du num\'ero \ref{subsection4.3}~;
au-dessus de la seconde couche de
$\mathcal C (C_k \ast C_\ell, \{\mu,\nu\})$,
on colle le long de chaque ar\^ete du premier prisme
un prisme infini sur un $\ell$-gone r\'egulier~;
etc.
\par

Nous laissons au lecteur le soin de se convaincre que
le diagramme ainsi obtenu est quasi-isom\'etrique
au produit de l'arbre biparti de degr\'es $k,\ell$
mentionn\'e en fin de \ref{subsection4.3} et d'une droite.
De plus, tous ces diagrammes 
(\`a la seule exception de celui correspondant \`a $k = \ell = 2$)
sont quasi-isom\'etriques entre eux.

\section*{4$^{bis}$. Invariants g\'eom\'etriques de Hopf \`a Gromov}

L'avantage de repr\'esenter un groupe $G$ par 
un diagramme de Cayley ou un graphe de Cayley,
c'est-\`a-dire par un objet g\'eom\'etrique, est de pouvoir d\'efinir pour $G$
des \emph{invariants g\'eom\'etriques}.
Ainsi, on peut voir le \emph{nombre de bouts} 
du graphe de Cayley d'un groupe de type fini
comme un invariant g\'eom\'etrique du groupe lui-m\^eme
(travaux de Hopf, et de ses \'etudiants Freudenthal et Specker 
dans les ann\'ees 1940).

Ces invariants g\'eom\'etriques jouent un r\^ole pr\'epond\'erant
en th\'eorie des groupes depuis les travaux de Gromov,
dont les tr\`es influents \cite{Grom--87} et \cite{Grom--93}.

\section{Sph\`eres d'homologie}
\label{section5}

Soit $K$ un n\oe{}ud, orient\'e, donn\'e par un  diagramme $\mathcal D$
dans le plan orient\'e.
Les points de croisement divisent  la projection de $K$ 
en un certain nombre d'arcs, chacun h\'eritant de l'orientation de $K$.
Notons $\sigma_1$ l'un d'entre eux (arbitraire),
et $\sigma_2, \hdots, \sigma_n$ les suivants dans l'ordre impos\'e par l'orientation~;
notons $s_1, \hdots, s_n$ les g\'en\'erateurs de Wirtinger correspondants de $G_K$.
Choisissons $s_1$ comme m\'eridien de $\partial E (K)$. 
Nous allons d\'ecrire une recette fournissant le parall\`ele correspondant.
\par

En parcourant $D$ \`a partir du milieu de $\sigma_1$,
on passe par $n$ points de croisement 
o\`u le n\oe{}ud passe en-dessous d'un arc.
Notons, dans l'ordre, $P_1, \hdots, P_n$ ces points~;
notons $\sigma_{i_j}$ l'arc qui passe par dessus en $P_j$~;
posons $\epsilon_j = 1$ [respectivement $\epsilon_j = -1$]
si $\sigma_{i_j}$ passe par dessus  de gauche \`a droite [resp. de droite \`a gauche].
Alors le parall\`ele cherch\'e s'\'ecrit
\begin{equation*}
p \, = \, s_{i_1}^{\epsilon_{i_1}} \cdots s_{i_n}^{\epsilon_{i_n}} (s_1)^k
\end{equation*}
o\`u $k$ est tel que $p$ soit homologiquement trivial dans $H_1(E(K))$,
c'est-\`a-dire tel que
$\epsilon_{i_1} + \cdots + \epsilon_{i_n} + k = 0$.

Pour une projection convenable du n\oe{}ud de tr\`efle $K$,
on trouve
\begin{equation*}
G_K \, = \, \langle a,b,c \hskip.1cm \vert \hskip.1cm ab = bc = ca \rangle
\tag{$P'_2$}
\end{equation*}
et 
\begin{equation*}
m \, = \, a , \hskip.5cm p \, = \, c^{-1}a^{-1}b^{-1}a^3 .
\end{equation*} 
Si $d = ab = bc = ca$ comme\footnote{A
l'\'echange pr\`es de $a$ et $b$.} 
plus haut, notons que
$b = a^{-1}d$, $c = d a^{-1}$, et $ada = abca = d^2$, de sorte que
\begin{equation*}
G_K \, = \, \langle a,d \hskip.1cm \vert \hskip.1cm ada = d^2 \rangle 
\tag{$P_5$}
\end{equation*}
et
\begin{equation*}
\aligned
p \, &= \, c^{-1}a^{-1}b^{-1}a^3 \, = \,
\left( ad^{-1} \right) a^{-1} \left( d^{-1} a \right) a^3 \, 
\\
&= \, ad^{-1} (ada)^{-1} a^5 \, = \, ad^{-3}a^5 \, = \, d^{-3}a^6 
\endaligned
\end{equation*}
(car $d^3$ est central dans $G_K$).

Soit $T'$ un tore solide.
Choisissons un m\'eridien $m'$ et un parall\`ele $p'$ sur $\partial T'$~;
notons que $m' = 1 \in \pi_1(T')$ et que $p'$ engendre $\pi_1(T')$.
Pour tout entier $k \in \mathbf Z$, 
soit $h_k : \partial T' \longrightarrow \partial V(K)$ un hom\'eomorphisme
tel que $h_k(m') = a p^{-k} \hskip.1cm [= mp^{-k}]$ et $h_k(p') = p$
(un tel hom\'eomorphisme est bien d\'efini \`a isotopie pr\`es).
Notons
\begin{equation*}
M_k \, = \, E(K) \cup_{h_k} T' 
\end{equation*}
le r\'esultat du recollement de $E(K)$ et $T'$ selon $h_k$~;
c'est une $3$-vari\'et\'e compacte connexe orientable.

Comme $d^3$ est central,
$ap^{-k}$ s'\'ecrit aussi $d^{3k}a^{1-6k}$.
Le th\'eor\`eme de Seifert -- van Kampen\footnote{L'article
de Seifert date de 1931 et celui de van Kampen de 1933.
Ce dernier \'etait une commande de Zariski,
qui voulait une m\'ethode simple pour analyser le groupe fondamental
du compl\'ementaire d'une courbe alg\'ebrique
dans le plan projectif complexe.
Pour la d\'etermination de groupes fondamentaux avant
ce th\'eor\`eme ``de Seifert -- van Kampen'', et notamment chez Poincar\'e,
voir le chapitre 3 de \cite{Gord--99} et
le n$^{\text{o}}$ 4.1 de \cite{Stil--93}.}
montre que le groupe fondamental de $M_k$ 
est isomorphe au produit libre
de $G_K$ et $\pi_1(T') = \langle p' \rangle$
avec amalgamation donn\'ee par $h_k$, c'est-\`a-dire que
\begin{equation*}
\aligned
\pi_1(M_k) \, &= \,
\langle a,d,p' \hskip.1cm \vert \hskip.1cm
ada = d^2 , \hskip.1cm ap^{-k} = 1 , \hskip.1cm p = p' \rangle
\\
&= \, \langle a,d \hskip.1cm \vert \hskip.1cm 
ada = d^2, \hskip.1cm d^{3k} = a^{6k-1} \rangle .
\endaligned
\end{equation*} 
Notons que ce groupe est parfait~; 
en effet, si $\equiv$ d\'esigne une \'egalit\'e dans $\pi_1(M_k)$
modulo les commutateurs, nous avons d'abord
\begin{equation*}
d \, = \, d^2d^{-1} \, = \, adad^{-1} \, = \, a^2(a^{-1}dad^{-1}) \, \equiv \, a^2 ,
\end{equation*}
donc $d^{3k} \equiv a^{6k} = d^{3k}a$ et $a \equiv 1$,
et ensuite $d \equiv ada = d^2$, donc $d \equiv 1$.

On obtient le r\'esultat suivant.
Voir \cite{Dehn--10}, pages 117--121 de \cite{DeSt--87},
et les pages 67--74 de  \cite{Rham--69}~;
voir aussi  le num\'ero 8.4 de \cite{Stil--93}.

\medskip

\noindent\textbf{Th\'eor\`eme sur les sph\`eres d'homologie.}
\emph{
Pour tout entier $k \in \mathbf Z$, 
la vari\'et\'e $M_k$ est une sph\`ere d'homologie. De plus
\begin{itemize}
\item[-]
pour $k = 0$, c'est la sph\`ere usuelle~;
\item[-]
pour $k = 1$, le groupe $\pi_1(M_1)$ est le groupe binaire de l'icosa\`edre,
qui est un groupe parfait \`a $120$ \'el\'ements~;
\item[-]
pour $k \ne 0,1$, le groupe $\pi_1(M_k)$ est infini~;
\item[-]
pour $k$ et $k'$ distincts dans $\mathbf Z$,
les groupes $\pi_1(M_k)$ et $\pi_1(M_{k'})$ sont non isomorphes.
\end{itemize}
}

\medskip

\footnotesize

D\'ej\`a dans une courte note  \cite{Dehn--07} pr\'ec\'edant \cite{Dehn--10}, 
Dehn avait construit  des sph\`eres d'homologie, par une m\'ethode diff\'erente.
Etant donn\'e deux n\oe{}uds  $K$ et $K'$ dans $\mathbf S^3$, 
soit $M$ la vari\'et\'e obtenue 
\`a partir des ext\'erieurs $E(K)$ et $E(K')$ 
en recollant leurs bords, avec un m\'eridien $m$ de $\partial E (K)$
identifi\'e \`a un parall\`ele $p'$ de $\partial E (K')$,
et de m\^eme $p$ identifi\'e \`a $m'$.
Alors $M$ est une sph\`ere d'homologie~;
de plus, si $K$ et $K'$ sont non triviaux, 
$\mathbf Z^2$ s'injecte dans le produit amalgam\'e
$\pi_1(M) = \pi_1(K) \ast_{\pi_1(\mathbf T^2)} \pi_1(K')$,
qui est donc un groupe infini.

\normalsize

\medskip

Le type de construction consistant \`a \^oter 
certains parties toriques d'une vari\'et\'e 
pour les recoller diff\'eremment
est la premi\`ere apparition de \emph{chirurgies} en topologie.
Citons deux expositions r\'ecentes sur l'usage de chirurgies
pour l'\'etudes des n\oe udes et des vari\'et\'es de dimension trois~:
\cite{Boye--02} et \cite{Gord--09}.
\par
Dans l'\'etude des vari\'et\'es de grandes dimensions,
la chirurgie est un processus de transformation contr\^ol\'ee des vari\'et\'es,
selon le sch\'ema
\begin{equation*}
M_1 =
M_0 \cup_{\mathbf S^p \times \mathbf S^q} (\mathbf S^p \times \mathbf D^{q+1})
 \, \Longrightarrow \, 
M_2 =
M_0 \cup_{\mathbf S^p \times \mathbf S^q} (\mathbf D^{p+1} \times \mathbf S^q)
\end{equation*}
(avec $p+q+1 = \dim M_0$).
C'est un ingr\'edient essentiel dans les travaux des ann\'ees 1960 et suivantes
sur les vari\'et\'es de grandes dimensions 
(Kervaire, Milnor, Browder, Novikov, Sullivan, Wall, ...).
Dans le cas particulier o\`u $p=q=1$, on retrouve
certaines\footnote{Si $(p,m)$ est un couple (parall\`ele, m\'eridien)
sur $\partial E (K) \approx  \partial(\mathbf S^1 \times \mathbf D^2)$
et $(m',p')$ un couple (m\'eridien, parall\`ele) sur le bord d'un tore solide
$T' = \mathbf D^2 \times \mathbf S^1$,
ce sont les chirurgies de Dehn associ\'ees \`a un hom\'eomorphisme
$h : \partial E (K) \longrightarrow \partial T'$ tel que
$h(pm^k) = m'$ et $h(m) = p'$, pour un entier $k \in \mathbf Z$.
La donn\'ee du n\oe{}ud $K$ et de l'entier $k$ correspond \`a celle d'un
\emph{n\oe{}ud enrubann\'e} (= framed knot).} 
chirurgies de Dehn
\begin{equation*}
\mathbf S^3 =
E(K) \cup_{\mathbf T^2} (\mathbf S^1 \times \mathbf D^2)
 \, \Longrightarrow \, 
E(K) \cup_{\mathbf T^2} (\mathbf D^2 \times \mathbf S^1) .
\end{equation*}
Mais le sujet des chirurgies de Dehn, propre \`a la dimension trois,
est bien diff\'erent de celui des chirurgies utilis\'es dans l'\'etude des vari\'et\'es de grandes dimensions, et le premier sujet
ne semble gu\`ere avoir servi d'inspiration au second.

\section*{5$^{bis}$. A propos de sph\`eres d'homologie}

Les sph\`eres d'homologie jouent un r\^ole important en topologie~;
voir  \cite{Save--02}.
Evoquons quelques points de l'histoire post\'erieure \`a Dehn, comme suit.
\par

Toute $3$-vari\'et\'e orient\'ee $M$ est le bord d'une $4$-vari\'et\'e $Q$,
comme l'ont not\'e ind\'ependamment Thom
\cite{Thom--51} et Rochlin \cite{Roch--51}~;
on peut supposer de plus que  $Q$ est simplement connexe et 
poss\`ede une structure spin.
Le groupe $H_2(Q, \mathbf Z)$ est sans torsion (parce que $\pi_1(Q) = 1$)~;
la forme d'intersection
\begin{equation*}
H_2(Q, \mathbf Z) \times H_2(Q, \mathbf Z) \longrightarrow \mathbf Z 
\end{equation*}
est unimodulaire (parce que $M$ est une sph\`ere d'homologie)
et paire (parce que $Q$ a une structure spin).
Cette forme a une signature qui est un multiple de $8$
(c'est une propri\'et\'e arithm\'etique des formes unimodulaires paires),
dont la classe modulo $16$ ne d\'epend que de $M$
et pas du choix de $Q$
(c'est un r\'esultat topologique profond de Rochlin).
On obtient donc un \emph{invariant de Rochlin
$\mu(M) \in \mathbf Z / 2 \mathbf Z$
des sph\`eres d'homologie enti\`ere}
(qui n'est qu'une toute petite partie d'une histoire aux multiples facettes
d\'evoil\'ee par Rochlin, voir \cite{GuMa--86}).
\par

Cet invariant s'est r\'ev\'el\'e plus tard \^etre la r\'eduction modulo $2$
d'un \emph{invariant de Casson}, \`a valeurs dans $\mathbf Z$.

\medskip

Si un groupe $G$ est le groupe fondamental
d'une $n$-sph\`ere d'homologie, avec $n \ge 3$,
alors $G$ est de pr\'esentation finie, parfait, et $H_2(G;\mathbf Z) = 0$.
Lorsque $n \ge 5$, r\'eciproquement, tout groupe ayant ces trois propri\'et\'es
est groupe fondamental d'une $n$-sph\`ere d'homologie \cite{Kerv--69}.
Lorsque $n  = 4$ ou $n = 3$, on ne conna\^{\i}t pas
de conditions n\'ecessaires et suffisantes.
\par

Citons encore un autre r\'esultat de Kervaire~:
\begin{itemize}
\item[] 
\emph{le seul groupe fini non trivial qui est groupe fondamental 
d'une $3$-sph\`ere d'homologie est le groupe binaire de l'icosa\`edre,
\`a $120$ \'el\'ements~;}
\end{itemize}
ainsi qu'une cons\'equence des travaux de Perelman~:
\begin{itemize}
\item[] 
\emph{une $3$-sph\`ere d'homologie \`a groupe fondamental fini non trivial
est diff\'eomorphe au quotient de $SU(2)$ par le groupe binaire de l'icosa\`edre.}
\end{itemize}

\section{Algorithme(s) de Dehn}
\label{section6}

Soit $\langle S \hskip.1cm \vert \hskip.1cm R \rangle$ 
une pr\'esentation finie d'un groupe $G$.
Notons $R_*$ l'ensem\-ble des relations obtenues \`a partir de celles de $R$
par permutations cycliques et inversion.
Soit $w$ un mot en les lettres de $S \cup S^{-1}$~;
notons $\vert w \vert$ la \emph{longueur} de $w$,
c'est-\`a-dire le nombre de ses lettres.

Une pr\'esentation finie $\langle S \hskip.1cm \vert \hskip.1cm R \rangle$
d'un groupe $G$ est une \emph{pr\'esentation de Dehn} 
si elle satisfait la condition suivante~: 
toute relation $r=1$ dans $R_*$
peut s'\'ecrire sous la forme $r \equiv u_rv_r^{-1}$, c'est-\`a-dire $u_r = v_r$,
o\`u $u_r$ et $v_r$ sont des mots en les $s_j$ et leurs inverses 
tels que $\vert v_r \vert < \vert u_r \vert$.

C'est par exemple le cas de la pr\'esentation usuelle du groupe fondamental
d'une surface close orientable de genre 2
\begin{equation*}
\pi_1(\Sigma_2) \, = \, 
\langle s_1, \hdots, s_4 \hskip.1cm \vert \hskip.1cm 
s_1 s_2 s_1^{-1}s_2^{-1}  s_3 s_4s_3^{-1} s_4^{-1}= 1 \rangle .
\end{equation*}
Pour cette pr\'esentation, $\vert R \vert = 1$ et 
$\vert R_* \vert = 16$. A titre d'exemple, 
$s_2s_1^{-1}s_2^{-1}s_3s_4 = s_1^{-1}s_4s_3$
est l'une des $16$ relations de $R_*$,
\'ecrite sous la forme $u_r = v_r$.
\par

Pour une pr\'esentation de Dehn, 
\emph{l'algorithme de Dehn} pour le probl\`eme du mot fonctionne comme suit~:
\'etant donn\'e un mot $w$ en les $s_j$, 
s'il contient une paire de lettres cons\'ecutives
de la forme $ss^{-1}$, on la supprime,
et s'il contient un $u_r$,
on obtient un nouveau mot plus court en rempla\c cant $u_r$ par $v_r$~;
on r\'ep\`ete l'op\'eration tant qu'elle est possible
(en choisissant arbitrairement la paire $ss^{-1}$ \`a supprimer
ou le $u_r$ \`a remplacer s'il y a plusieurs possibilit\'es)~;
le mot $w$ repr\'esente $1$ ou non selon que
le proc\'ed\'e s'ar\^ete avec le mot vide ou avec un mot non vide.

Dehn a montr\'e que les pr\'esentations usuelles 
des groupes fondamentaux de surfaces sont des
``pr\'esentations de Dehn'' de sorte que le probl\`eme du mot
est r\'esoluble pour ces groupes.
De plus, dans un tel groupe, chaque \'el\'ement donn\'e
comme produit des g\'en\'erateurs peut \^etre algorithmiquement
transform\'e en une ``forme normale'', de sorte que
le probl\`eme de conjugaison est aussi r\'esoluble.
Dans un premier article \cite{Dehn--11},
Dehn \'etablit ces r\'esultats en utilisant les propri\'et\'es g\'eom\'etriques
(courbure n\'egative) du plan hyperbolique,
qui est le rev\^etement universel d'une surface de genre $g \ge 2$~;
dans un second article \cite{Dehn--12},
il utilise de mani\`ere essentielle la combinatoire des chemins ferm\'es
dans le diagramme de Cayley.
Dehn insiste sur la signification g\'eom\'etrique de ses r\'esultats~;
voir (WP$_{\text{top}}$) et (CP$_{\text{top}}$),
ci-dessus au chapitre \ref{section2}.
\par

Ces articles ont marqu\'e toute la suite du sujet,
dont les m\'ethodes diagrammatiques de van Kampen (1933),
pass\'ees inaper\c cues \`a l'\'epoque mais red\'ecouvertes par Lyndon (1966),
et par l\`a toute la th\'eorie de la petite simplification.
La terminologie ``algorithme de Dehn'', due \`a Magnus, 
appara\^{\i}t dans la th\`ese de son \'etudiant Martin Greendlinger (1960),
o\`u il est montr\'e que l'algorithme fonctionne pour des pr\'esenta\-tions 
satisfaisant \`a des conditions convenables de petite simplification
(selon la page 21 de \cite{ChMa--82}).
\par

Gromov a montr\'e qu'un groupe poss\`ede
une pr\'esentation de Dehn si et seulement s'il est hyperbolique,
au sens de \cite{Grom--87}~;
voir par exemple \cite{BrHa--99}, th\'eor\`eme III.$\Gamma$.2.6, page 450
(o\`u il conviendrait d'introduire $R_*$).

\section*{6$^{bis}$. Fonctions de Dehn}

Soit $\langle S \hskip.1cm \vert \hskip.1cm R \rangle$ 
une pr\'esentation finie d'un groupe $G$.
Soit $w$ un mot en les lettres de $S \cup S^{-1}$
qui repr\'esente l'\'el\'ement neutre de $G$,
ce qu'on \'ecrit $w =_G 1$~;
il existe donc des relations $r_i \in R_*$
et des mots $v_i$, tels que $w = \prod_{i=1}^N v_i r_i v_i^{-1}$.
Par d\'efinition, \emph{l'aire} de $w$ est le minimum $A(w)$
des entiers $N$ pour lesquels une telle \'ecriture est possible. 
La \emph{fonction de Dehn} $\delta : \mathbf N \longrightarrow \mathbf N$
est d\'efinie par
\begin{equation*}
\delta(n) \, = \, \max \big\{ A(w) \hskip.1cm \big\vert \hskip.1cm
w =_G 1 \hskip.2cm \text{et} \hskip.2cm \vert w \vert \le n \big\} .
\end{equation*}
A priori, $\delta$ d\'epend de la pr\'esentation donn\'ee~;
toutefois, on montre facilement que, si deux pr\'esentations finies
d'un m\^eme groupe donnent lieu \`a des fonctions $\delta_1$ et $\delta_2$,
il existe une constante $C > 0$ telle que 
\begin{equation*}
\delta_i(n) \, \le \, C \delta_j ( Cn+C) + Cn + C
\hskip.5cm \text{pour tout} \hskip.2cm n \ge 0
\end{equation*}
($\{i,j\} = \{1,2\}$).
Dans ce sens, la fonction de Dehn
ne d\'epend \`a \'equivalence pr\`es que du groupe de pr\'esentation finie $G$,
et pas de la pr\'esentation finie choisie.
Les fonctions de Dehn permettent de quantifier 
le degr\'e de solubilit\'e du probl\`eme du mot pour un groupe de pr\'esentation finie~:
\begin{itemize}
\item[-]
un groupe a un probl\`eme du mot r\'esoluble si et seulement si
sa fonction de Dehn est r\'ecursive~;
\item[-]
un groupe a une fonction de Dehn lin\'eaire si et seulement s'il est
hyperbolique au sens de Gromov \cite{Grom--87}~;
\item[-]
si $IP$ d\'esigne 
l'ensemble des nombre r\'eels  $\rho \ge 1$
tels qu'il existe un groupe de pr\'esentation finie dont la fonction de Dehn
est \'equivalente \`a $n \longmapsto n^{\rho}$
(avec ``IP'' pour ``Isoperimetric Spectrum''),
alors l'adh\'erence de $IP$ est \'egale \`a $\{1\} \cup [2,\infty[$~;
\item[-] il existe un groupe \`a un relateur, 
donc \`a  probl\`eme du mot r\'esoluble (c'est un r\'esultat de Magnus), 
dont la fonction de Dehn est donn\'ee \`a \'equivalence pr\`es par
$\delta(n) = \epsilon_n(n)$, o\`u les fonctions $\epsilon_n$
sont d\'efinies r\'ecursivement par
\begin{equation*}
\epsilon_0(k) = k , \hskip.3cm
\epsilon_1(k) = 2^k , \hskip.3cm
\hdots , \hskip.3cm
\epsilon_n(k) = 2^{\epsilon_{n-1}(k)} , \hskip.3cm
\hdots ,
\end{equation*}
de sorte que $\delta$ cro\^{\i}t plus vite que toute exponentielle it\'er\'ee~;
\end{itemize}
pour tous ces r\'esultats, voir l'exposition de \cite{Brid--02}.
On comprend aussi presque quelles sont les fonctions 
$\mathbf N \longrightarrow \mathbf N$ qui sont (\`a \'equivalence pr\`es)
des fonctions de Dehn~: voir \cite{SaBR--02} et \cite{BORS--02}.

\section{Groupes de pr\'esentation finie, $2$-complexes,
\\ 
$3$- et $4$-vari\'et\'es}
\label{section7}

Dehn a vu que  
\begin{itemize}
\item[(*)]
\emph{tout groupe de pr\'esentation finie 
est groupe fondamental d'un $2$-complexe, 
et aussi d'une vari\'et\'e compacte orientable de dimension quatre~;}
\end{itemize}
voir le \S~2 du chapitre I de \cite{Dehn--10} 
et le d\'ebut du chapitre III de \cite{Dehn--11}.
Voici, en termes d'\emph{aujourd'hui},
les pas d'une d\'emonstration de (*)
inspir\'ee des articles de Dehn.

\medskip

(i) A une pr\'esentation finie d'un groupe $G$,
on associe un CW-complexe fini $C$, 
qui  a une seule $0$-cellule,
des $1$-cellules  en bijection avec les g\'en\'erateurs de la pr\'esentation,
et des $2$-cellules en bijection avec les relations, 
attach\'ees au $1$-squelette conform\'ement aux relations~;
il en r\'esulte que $\pi_1(C) = G$.
Le rev\^etement universel de $C$ est le \emph{complexe de Cayley}
de la pr\'esentation donn\'ee, et son $1$-squelette 
s'identifie au graphe de Cayley de la pr\'esentation.
Une subdivision ad hoc permet de consid\'erer $C$
comme un complexe simplicial fini.
\par

(ii) Un $2$-complexe simplicial fini se plonge dans $\mathbf R^5$,
ce dont on s'assure facilement 
par un argument de position g\'en\'erale~:
en dimension~$5$, 
deux plans affines peuvent toujours \^etre rendus disjoints 
par un d\'eplacement arbitrairement petit de l'un d'eux.
\par

(iii) $C$ poss\`ede un \emph{voisinage ouvert r\'egulier} $V$
dans $\mathbf R^5$, dont $C$ est un r\'etracte par d\'eformation,
de sorte que $\pi_1(C) \approx \pi_1(V) \approx \pi_1(\overline{V})$.
Comme on peut munir $\mathbf R^5$ d'une structure 
de complexe simplicial dont $C$ soit un sous-complexe,
ceci est un cas particulier du r\'esultat g\'en\'eral \'el\'ementaire suivant~:
si $L$ est un sous-complexe d'un complexe simplicial $K$,
et si $V(L)$ d\'esigne la r\'eunion des int\'erieurs des simplexes
de $K$ contenant un sommet de $L$ dans une subdivision barycentrique de $L$,
alors $V(L)$ est un voisinage de $K$ 
qui se r\'etracte par d\'eformation 
sur $L$~;
voir le no 9 du chapitre II dans \cite{EiSt--52}.
\par

(iv) D'une part, comme $C$ est de codimension trois dans $\overline{V}$,
nous avons par position g\'en\'erale un isomorphisme 
$\pi_1(\overline{V}) \approx \pi_1(\overline{V} \smallsetminus C)$.
D'autre part, $\overline{V} \smallsetminus C$ se r\'etracte par d\'eformation
sur le bord $\partial \overline{V}$ de $\overline{V}$, qui est une sous-vari\'et\'e close
de codimension $1$ dans $\mathbf R^5$, en particulier
une vari\'et\'e orientable de dimension $4$.
On a donc en fin de compte 
\begin{equation*}
G \, \approx \, \pi_1(C) \, \approx \, \pi_1(\partial V) .
\end{equation*}

\medskip

\footnotesize

\subsection{Variante}

A la place du fait invoqu\'e au point (ii) ci-dessus,
on peut utiliser le fait suivant,
cas  particulier d'un \'enonc\'e du d\'ebut de \cite{Stal--65}~:
tout $2$-complexe fini dont le $1$-squelette est planaire,
en particulier tout $2$-complexe fini ayant une seule $0$-cellule
se plonge dans $\mathbf R^4$.
En adaptant (iii) et (iv), 
on obtient un voisinage ouvert r\'egulier $V'$ dans $\mathbf R^4$
tel que le groupe fondamental $\pi_1(\partial V')$ 
se surjecte naturellement sur $\pi_1(V') \approx \pi_1(C)$.
Notons $W = V' \sqcup_{\partial V'} V'$ le double de $V'$,
c'est-\`a-dire la vari\'et\'e close de dimension $4$
obtenue en recollant deux copies de $V'$ sur leur bord commun~;
le th\'eor\`eme de Seifert -- van Kampen implique
que $\pi_1(W)$ est isomorphe \`a la somme amalgam\'ee
$\pi_1(V') \ast_{\pi_1(\partial V')} \pi_1(V')$,
c'est-\`a-dire \`a $\pi_1(V') \approx \pi_1(C) \approx G$.
Je remercie Martin Bridson pour m'avoir signal\'e \cite{Stal--65}.

\subsection{Digression}

Consid\'erons l'exemple du groupe cyclique d'ordre $7$
donn\'e par la pr\'esentation 
\begin{equation*}
\langle a,b,c \hskip.1cm \vert \hskip.1cm
a^2b^{-1} = bac^{-1} = cac = c^{-1}ab = b^2c = 1 \rangle
\end{equation*}
(dont on d\'eduit facilement $b = a^2$, $c = a^3$, 
de sorte que ce groupe, qui a aussi la pr\'esentation
$\langle a \hskip.1cm \vert \hskip.1cm a^7 = 1 \rangle$,
est bien cyclique d'ordre $7$).
Le complexe $C$ correspondant 
se plonge comme ci-dessus dans $\mathbf R^4$.
On peut v\'erifier que son rev\^etement universel $\widetilde C$
est le $2$-squelette du $6$-simplexe,
qui est l'un des exemples standards d'un $2$-complexe
\emph{qui ne se plonge pas} dans $\mathbf R^4$
(ce qui est connu depuis les travaux de van Kampen et Flores des ann\'ees 1930,
voir \cite{FrKT--94}).
Nous avons ainsi obtenu un exemple
de $2$-sous-complexe fini de $\mathbf R^4$
dont le rev\^etement universel ne se plonge pas dans $\mathbf R^4$,
ce qui (\`a tort ou \`a raison) surprend l'auteur.
\par

\subsection{Autre argument pour (*)}

Soit $M$ la somme connexe  de $n$ copies
de $\mathbf S^1 \times \mathbf S^3$~; 
son groupe fondamental  est libre
sur $n$ g\'en\'erateurs $s_1, \hdots, s_n$.
Etant donn\'e $m$ relations $r_1, \hdots, r_m$ en les $s_j$,
on repr\'esente chaque $r_i$ par un lacet dans $M$, 
ayant un voisinage tubulaire $C_i$ hom\'eomorphe 
\`a $\mathbf S^1 \times \mathbf D^3$~;
noter que le bord de $\mathbf S^1 \times \mathbf D^3$
est hom\'eomorphe au bord de la vari\'et\'e simplement connexe
et $\mathbf D^2 \times \mathbf S^2$.
On modifie alors $M$ ``par chirurgie'',
en rempla\c cant chaque $C_i$
par une copie de $\mathbf D^2 \times \mathbf S^2$.
Le groupe fondamental de la vari\'et\'e ainsi obtenue
est isomorphe au groupe de pr\'esentation
$\langle s_1, \hdots, s_n \hskip.1cm \vert \hskip.1cm r_1, \hdots, r_m \rangle$.
Voir la page 187 de \cite{SeTh--80}.
 
\normalsize

\subsection{Dimensions $\le 3$ et $\ge 4$}

Dehn ayant montr\'e que tout groupe de pr\'esentation finie
est groupe fondamental d'une $4$-vari\'et\'e close,
on peut se demander pourquoi la th\'eorie combinatoire des groupes,
qui est ins\'eparable d\`es ses premiers d\'ebuts de la topologie de dimension trois,
a pendant longtemps bien moins interagi
avec l'\'etude des vari\'et\'es de dimension $\ge 4$
(exemples~: les surfaces complexes~!).
Faute de pouvoir discuter s\'erieusement la question,
notons que  ce n'est que dans les ann\'ees 1920 que parut
le livre de Lefschetz sur
\emph{l'analysis situs et la g\'eom\'etrie alg\'ebrique} \cite{Lefs--24}.
Lefschetz qui \'ecrira plus tard que
ce fut ``son lot de planter le harpon de la topologie alg\'ebrique
dans le corps de la baleine de la g\'eom\'etrie alg\'ebrique''
(page 854 de \cite{Lefs--68}).

\section{Allusion \`a d'autres contributions}
\label{section8}

Nous avons d\'ej\`a fait allusion \`a
\begin{itemize}
\item[(i)]
la solution de Dehn du 3e probl\`eme de Hilbert
(au chapitre~\ref{section1})~;
\item[(ii)]
sa d\'emonstration du fait 
qu'il n'existe pas de d\'eformation continue
de $\mathbf R^3$ transformant un n\oe{}ud de tr\`efle gauche 
en un n\oe{}ud de tr\`efle droit, en d'autres termes du fait
que le n\oe{}ud de tr\`efle est ``chiral''
(au chapitre \ref{section2}).
\end{itemize}
Sans entrer dans aucun d\'etail, mentionnons tout de m\^eme encore~:
\begin{itemize}
\item[(iii)]
le th\'eor\`eme de Dehn-Nielsen concernant une surface close $S$, 
selon lequel tout automorphisme ext\'erieur de $\pi_1(S)$
est induit par un hom\'eo\-morphisme de $S$ (vers 1921)~;
\item[(iv)]
le r\^ole des ``twists de Dehn'' dans l'\'etude \cite{Dehn--38}
des groupes de diff\'eo\-morphismes 
et des groupes modulaires de surfaces (= mapping class groups)~;
\item[(v)]
et les \'equations de Dehn-Sommerville, portant sur les nombres
des faces de dimensions $0, 1, \hdots d$
d'une triangulation de la sph\`ere de dimension $d$,
\'equations conjectur\'ees par Dehn (1905)
et d\'emontr\'ees par Sommerville (1927).
\end{itemize}

\section{Quelques souvenirs et commentaires}
\label{section9}

\begin{center}
De Carl Ludwig Siegel \cite{Sieg--64}
\end{center}

``J'appris (...) quelle rare chance c'est que d'avoir des coll\`egues acad\'emiques
qui travaillent  solidairement, sans \'ego\"{\i}sme et
sans l'id\'ee d'une ambition personnelle, 
plut\^ot que d'\'emettre des directives venant de leurs positions hautaines.''
\par

``Il avait un esprit philosophique au sens de Schiller et,
comme il \'etait \'epris de contradiction,
une conversation avec lui devenait souvent  une discussion fructueuse.
Il s'int\'eressait beaucoup \`a l'histoire, ancienne et moderne (...).''
\par

A propos du \emph{Freiheitssatz}, 
propos\'e par Dehn comme sujet de th\`ese \`a Magnus, 
Siegel \'ecrit aussi~:
``qu'il en d\'ecouvrit une d\'emonstration, et \`a l'occasion
disait \`a ses amis comment elle allait''.
\footnote{
Elle \'etait  g\'eom\'etrique,
contrairement \`a celle de Magnus,
dont le caract\`ere alg\'ebrique
valut \`a Magnus que Dehn lui serve  le commentaire suivant~:
``Da sind Sie also blind gegangen'', ou
``Ainsi avez-vous proc\'ed\'e \`a l'aveugle''
(\cite{Magn--78}, page 139).
}

\smallskip
\begin{center}
Du cin\'easte Arthur Penn  \cite{Penn}
\end{center}

\noindent
Penn \'etait \`a Black Mountain en 1948.

``But then world events occurred, of a certain extraordinary character, 
which was that in Germany there was the rise of Hitler, 
and the creation of an enormous group of refugees 
who were some of the most remarkably talented people in all of Europe
who were suddenly persona non grata, and forced out. 
And, although the United States had something of a welcoming policy, 
it was a peculiar period of, yeah we'll let you come into the United States, 
but we won't accept any of your credentials from Europe, 
because somehow that was just not the American way, 
so an awful lot of doctors for instance came in and had to take their medical exams over again. Psychoanalysts of world wide reputation had to go and sit for another medical exam. 
And we at Black Mountain were very fortunate in that people like Gropius and Albers, 
and Max Dehn who was a marvelous physicist and mathematician, came. 
So that a large part of our faculty population was made up 
of these absolutely extraordinary people who were able to make it to America, 
but couldn't get a job at an American university. 
Although they'd been pillars of the Bauhaus in Germany, 
which was a remarkable place.''

\smallskip
\begin{center}
Dernier alin\'ea de la pr\'eface du livre \cite{MaKS--66} 
\end{center}

``This book is dedicated to the memory of Max Dehn. 
We believe this to be more than an acknowledgment of a personal indebtedness 
by one of the authors who was Dehn's student. 
The stimulating effect of Dehn's ideas on presentation theory 
was propagated not only through his publications, 
but also through talks and personal contacts; 
it has been much greater than can be documented by his papers. 
Dehn pointed out the importance of fully invariant subgroups in 1923 in a talk 
(which was mimeographed and widely circulated but never published). 
His insistence on the importance of the word problem, 
which he formulated more than fifty years ago, 
has by now been vindicated beyond all expectations.''

\smallskip 
\begin{center}
De Wilhelm Magnus \cite{Magn--78} 
\end{center}

``We can also say that being a mathematician 
was an essential part of his personality 
and that it influenced also his very well founded and deep interests 
in the humanities, in art, and in nature.''

\smallskip
\begin{center}
De John Stillwell \cite{DeSt--87}, page 253
\end{center}

\noindent
Extrait de l'introduction \`a la traduction de \cite{Dehn--38},
article d\'ej\`a mentionn\'e au chapitre \ref{section8},
qui reprenait et d\'eveloppait un expos\'e de 1922
sur les groupes modulaires de surfaces.
\smallskip

``Perhaps because of its very demanding proof,
the result went unnoticed until it was rediscovered independently by Lickorish (1962).''

``Dehn had another idea which also went unnoticed, 
even when it was rediscovered in the 
unpublished\footnote{Publi\'e une ann\'ee 
apr\`es la parution de \cite{DeSt--87}.},  but famous, 
Thurston \cite{Thur--76}.
This was the idea of studying the mapping class group
 by its action on the space of simple curve systems.''

\smallskip\begin{center}
D'Andr\'e Weil \cite{Weil--91}, pages 52--53
\end{center}

``J'ai rencontr\'e deux hommes dans ma vie 
dont le souvenir me fait penser \`a Socrate~;
ce sont Max Dehn et Brice Parain\footnote{Philosophe et essayiste fran\c cais,
1897--1971~; ami d'Albert Camus~; appara\^{\i}t  dans
``Vivre sa vie'' de Jean-Luc Godard \cite{Goda--62}.
Puissent deux des phrases qu'il y prononce 
s'appliquer au pr\'esent texte~:
``Le mensonge c'est un des moyens de la recherche'' et
``Il faut passer par l'erreur pour arriver \`a la v\'erit\'e''.}.
Ils avaient de Socrate, tel que nous l'imaginons 
d'apr\`es le t\'emoignage de ses disciples,
le rayonnement qui fait qu'on s'incline naturellement devant leur m\'emoire~; 
c'est l\`a une qualit\'e \`a la fois intellectuelle et morale 
que le mot de $<<$sagesse$>>$ est peut-\^etre le mieux fait pour exprimer, 
car la saintet\'e est autre chose. 
A c\^ot\'e du sage, le saint n'est peut-\^etre qu'un sp\'ecialiste --
un sp\'ecialiste de la saintet\'e~; le sage
n'a pas de sp\'ecialit\'e.
Ceci n'est pas \`a dire, loin de l\`a, 
que Dehn n'ait pas \'et\'e un math\'ematicien de grand talent~;
il a laiss\'e une \oe{}uvre de haute qualit\'e. 
Mais pour un tel homme la v\'erit\'e est une, 
et la math\'ematique n'est que l'un des miroirs o\`u elle se r\'efl\`ete, 
avec plus de puret\'e peut-\^etre que dans d'autres. 
Esprit universel, il avait une profonde connaissance 
de la philosophie et de la math\'ematique grecques.''

\section{Autour de 1910 -- floril\`ege}
\label{section10}

1907  \emph{Les Demoiselles d'Avignon,}
de Picasso~; titre original~: El Burdel de Avi\~n\' on.

1907-1909 et 1916
Cours \`a Gen\`eve par Ferdinand de Saussure (1857-1913)
et publication par ses \'el\`eves du
\emph{Cours de linguistique g\'en\'erale.}

30 juin 1908 Chute d'une m\'et\'eorite g\'eante \`a Toungouska, en Sib\'erie.

1908 Quatri\`eme congr\`es international des math\'ematiciens, Rome.
Dehn et sa femme y participent.

1908--1914 Andr\'e Gide cr\'eateur et premier directeur de la
\emph{Nouvelle Revue fran\c caise.}

1909 Publication par Walter Ritz (1878-1909) 
de travaux de physique math\'ematique,
dans lesquels il met au point une m\'ethode de calcul,
rapidement reprise en Russie par Galerkin, 
qui deviendra c\'el\`ebre sous le nom de
``m\'ethode de Ritz-Galerkin'',
et qui est \`a la base de 50 \% 
de tous les calculs scientifiques d'aujourd'hui \cite{GaWa}.

1909 En juillet, un d\'ecret du cardinal Respighi, vicaire 
de Rome,
interdit aux eccl\'esiastiques de la ville sainte
la fr\'equentation du cin\'ema\-tographe, sous peine de suspension
\emph{a divinis}.
(1913 Ouverture au Vatican du cin\'ema pontifical r\'eserv\'e
aux membres du clerg\'e.)

1909 Premi\`ere travers\'ee de la Manche en avion, le 25 juillet, par Louis Bl\'eriot.

Vers 1910 Travaux de Luitzen Egbertus Jan Brouwer en topologie, dont~:
\emph{Beweis der Invarianz der Dimensionenzahl,} 
Math. Ann. \textbf{70} (1911), 161--165~;
\emph{\"Uber Abbildung von Mannigfaltigkeiten,}
Math. Ann. \textbf{71} (1912), 97--115
(degr\'e de Brouwer, th\'eor\`eme du point fixe).

8 mars 1910 Cr\'eation \`a Copenhague d'une ``Journ\'ee de la femme''
par une conf\'ed\'eration internationale de femmes socialistes,
en vue de faire admettre le vote des femmes.
Vote  acquis en 1893 en Nouvelle-Z\'elande, 
en 1918 en Allemagne, 
en 1944-45 en France, et en
1971 en Suisse (ce n'est m\^eme pas un record). 

1910 \emph{Premier quatuor \`a cordes}
en la mineur de B\'ela Bartok, cr\'e\'e le 19 mars \`a Budapest.  

1910  \emph{Khodynka}, dernier \'ecrit litt\'eraire de L\'eon Tolsto\"{\i}, 
mort le 7 novembre.

1910 \emph{Aquarelle abstraite} de Vassily Kandinsky.

1910 Th\`eses de Erich Hecke et Richard Courant.
64 \'etudiants de Hilbert  obtinrent leurs th\`eses entre 1898 et 1915,
dont Dehn en 1900~; et 10 autres apr\`es l'interruption de la guerre
(d'apr\`es le ``Mathematics Genealogy Project'').
Pour comparaison~: Selim Krein, 81 \'etudiants de th\`ese~;
Andrei Kolmogorov, 79 \'etudiants~;
Charles Ehresmann, 76 \'etudiants~;
Wilhelm Magnus, 73 \'etudiants~; 
Beno Eckmann, 62 \'etudiants~;
Heinz Hopf, 49 \'etudiants~;
Henri Cartan, 14 \'etudiants~;
etc, voir  {\tt http://genealogy.math.ndsu.nodak.edu/extrema.php}

1910--1913 Cr\'eations parisiennes d'Igor Stravinski pour les ballets de Diaghilev~:
\emph{L'oiseau de feu}  le 25 juin 1910,
\emph{Petrouchka}  le 13 juin 1911 et 
le \emph{Sacre du Printemps}  le 29 mai 1913.

1910--1913 \emph{Principia Mathematica,} 
de Alfred North Whitehead et Bertrand Russel
(trois volumes).

1910 Fondation de la Soci\'et\'e Math\'ematique suisse,
apr\`es celles des
London Mathematical Society (1865), 
Soci\'et\'e Math\'ematique de France (1872),
New York Mathematical Society (1888)
devenue l'American Mathematical Society (1894), 
Deutsche Mathematiker-Vereinigung (1890),
et \"Osterreichische Mathematische Gesellschaft (1903), 
et avant celles des
Unione Mathematica Italiana (1922)
et European Mathematical Society (1990).

1911 Frictions franco-allemandes 
au Maroc, le ``coup d'Agadir'' le 1er juillet.

1911 Roald Amundsen atteint le p\^ole Sud, le 14 d\'ecembre, un mois
avant Robert Falcon Scott (qui mourra dans une temp\^ete au retour).

1912  Cr\'eation le 16 octobre du \emph{Pierrot lunaire}, 
Op.~18 d'Arnold Sch\"onberg
(\oe{}uvre atonale, annon\c cant le dod\'ecaphonisme).

1912 Naufrage du Titanic le 14 avril.

1912 Franz Kafka commence \`a r\'ediger \emph{la M\'etamorphose}.

1912 Cinqui\`eme congr\`es international des math\'ematiciens, Cambridge (G.B.).

1913 Publication de  \emph{Du c\^ot\'e de chez Swann} de Marcel Proust,
premi\`ere partie de  \emph{\`A~la recherche du temps perdu.}

1913 
Derni\`ere course de la diligence Aigle Les Diablerets, le 22 d\'e\-cembre et
mise en exploitation du chemin de fer \'electrique.

1913 On compte 1370 automobiles \`a Gen\`eve.

1913 \emph{Die Idee der Riemannschen Fl\"ache,} de Hermann Weyl.

1914 Cr\'eation du personnage de Charlot dans \emph{Charlot vagabond}.

1914 \emph{Les gens de Dublin} de James Joyce.

1914 \emph{Grundz\"uge der Mengenlehre,} de Felix Hausdorff.

\medskip

\begin{center}
Quelques empereurs, rois et pr\'esidents autour de 1910
\end{center}

\textsc{Guillaume II}, 
dernier empereur allemand 
et dernier roi de Prusse de 1888 \`a 1918.

\textsc{Edouard VII} 
roi d'Angleterre de 1901 \`a 1910~;
\textsc{Georges V}
de 1910 \`a 1936.

\textsc{Aixinjueluo Puyi}, douzi\`eme et dernier empereur de la dynastie Qing
de 1908 \`a 1912~;
\textsc{Sun Yat-sen} pr\'esident de la R\'epublique de Chine en 1912.

Theodore \textsc{Roosevelt}, pr\'esident r\'epublicain des Etats-Unis de 1901 \`a 1909~;
William Howard \textsc{Taft}, idem de 1909 \`a 1913~;
Thomas Woodrow \textsc{Wilson}, pr\'esident d\'emocrate de 1913 \`a 1921.

Armand \textsc{Falli\`ere} pr\'esident 
de la Troisi\`eme R\'epublique Fran\c caise de 1906 \`a 1913, 
Raymond \textsc{Poincar\'e}
(Raymond, 1860--1934,
cousin du math\'e\-maticien pr\'enomm\'e Henri, 1854-1912.)

\textsc{Victor-Emmanuel III} 
roi d'Italie de 1900 \`a 1944.

\textsc{Nicolas II} 
dernier tsar de toutes les Russies de 1894 \`a 1917.

Robert \textsc{Comtesse}, pr\'esident de la conf\'ed\'eration suisse en 1904 et 1910.

 \newpage

\begin{center}
Dates de quelques math\'ematiciens mentionn\'es dans l'article
\end{center}

\par\noindent
Arthur \textsc{Cayley}, 1821--1895~;
fils\footnote{Au sens du ``Mathematics Genealogy Project''.
Presque toutes nos dates sont tir\'ees du site correspondant~:
{\tt  http://genealogy.math.ndsu.nodak.edu/}
}
de Hopkins.

\par\noindent
Peter Guthrie \textsc{Tait}, 1831--1901.

\par\noindent
Georg Ferdinand Ludwig Philipp \textsc{Cantor}, 1845--1918~;
\par
fils de Kummer et Weierstrass.

\noindent
Felix \textsc{Klein}, 1849--1925~;
fils de Pl\"ucker et Lipschitz.

\par\noindent
Henri \textsc{Poincar\'e}, 1854-1912~;
fils de Hermite.

\par\noindent
Walther \textsc{von Dyck}, 1856--1934~; fils de Klein.

\par\noindent
Karl Emmanuel Robert \textsc{Fricke}, 1861--1930~; fils de Klein.

\par\noindent
David \textsc{Hilbert}, 1862--1943~; fils de Lindemann.

\par\noindent
Wilhelm \textsc{Wirtinger}, 1865--1945~;
fils de Weyr et von Escherich.

\par\noindent
Poul \textsc{Heegaard}, 1871--1948~; th\`ese\footnote{Remarquablement
sans patron de th\`ese, bien que Heegaard ait \'et\'e influenc\'e
par ses contacts directs avec \textsc{Klein} et Julius \textsc{Petersen} en 1998,
et qu'il ait lu (et corrig\'e~!) l'Analsis situs de \textsc{Poincar\'e}.}
en 1898.

\par\noindent
Walter  \textsc{Ritz}, 1878--1909~;
influenc\'e par Hilbert~;
\par th\`ese avec Woldemar Voigt (physicien).

\par\noindent
$\Rightarrow$ \textsc{Max Dehn, 1878--1952}~; fils de Hilbert.

\par\noindent
Heinrich Franz Friedrich \textsc{Tietze}, 1880--1964~; fils de von Escherich.

\par\noindent
Solomon \textsc{Lefschetz}, 1884--1972~; fils de Story.

\par\noindent
Ludwig \textsc{Bieberbach}, 1886--1982~; fils de Klein.

\par\noindent
James Waddell \textsc{Alexander} II, 1888--1971~; fils de Veblen.

\par\noindent
Jakob \textsc{Nielsen}, 1890--1959~; fils de  Landsbert et Dehn.

\par\noindent
Heinz \textsc{Hopf}, 1894--1971~;
fils de Erhard Schmidt  et Bieberbach.

\par\noindent
Carl Ludwig \textsc{Siegel}, 1896--1981~; fils de Landau.

\par\noindent
Emil Leon \textsc{Post}, 1897--1954~; fils de Keyser.

\par\noindent
Hellmuth \textsc{Kneser}, 1898--1973~; fils de Hilbert.

\par\noindent
Oscar Ascher \textsc{Zariski}, 1899--1986~; fils de Castelnuovo.

\par\noindent
Pyotr Sergeyevich \textsc{Novikov}, 1901--1975~; fils de Luzin.

\par\noindent
Andrei Andreyevich \textsc{Markov}, junior\footnote{Fils biologique de
Andrei Andreyevich Markov, 1856--1922~;
ce dernier fils math\'ematique de Chebychev.}, 
1903--1979.

\par\noindent
Georges \textsc{de Rham}, 1903--1990~;
fils de Lebesgue.

\par\noindent
Alonzo \textsc{Church}, 1903--1995~; fils de Veblen.

\par\noindent
Hans \textsc{Freudenthal}, 1905--1990~;
fils de Heinz Hopf.

\par\noindent
Ruth \textsc{Moufang}, 1905--1977~; fille de Dehn.

\par\noindent
Kurt \textsc{G\"odel}, 1906--1978~; fils de Hahn.

\par\noindent
Werner \textsc{Burau}, 1906--1994~;  fils de Reidemeister.

\par\noindent
Ott-Heinrich \textsc{Keller}, 1906--1990~; fils de Dehn.

\par\noindent
Andr\'e \textsc{Weil}, 1906--1998~; fils de Hadamard et Picard.

\par\noindent
Erich \textsc{K\"ahler}, 1906--2000~; fils de Lichtenstein.

\par\noindent
Wilhelm \textsc{Magnus}, 1907--1990~; fils de Dehn.

\par\noindent
Karl Johannes Herbert \textsc{Seifert}, 1907--1996~;
\par
fils de Threlfall et van der Waerden.

\par\noindent
Egbert Rudolf \textsc{van Kampen}, 1908--1942~;
fils de van der Woude.

\par\noindent
Alan \textsc{Turing}, 1912--1954~; fils de Church.

\par\noindent
Christos \textsc{Papakyriakopoulos}, 1914--1976.

\par\noindent
Roger Conant \textsc{Lyndon}, 1917--1988~;
fils de Mac Lane.

\par\noindent
Vladimir \textsc{Rochlin}, 1919--1984~;
fils de Kolmogorov et Pontryagin.

\par\noindent
Graham \textsc{Higman}, 1917--2008~;
fils de J.H.C. Whitehead.

\par\noindent
William Werner \textsc{Boone}, 1920--1983~; fils de Church.

\par\noindent
Ernst Paul \textsc{Specker}, 1920~; fils de Heinz Hopf et Beno Eckmann.

\par\noindent
Ren\'e \textsc{Thom}, 1923--2002~; fils de Henri Cartan.

\par\noindent
John Leslie \textsc{Britton}, 1927--1994~;
fils de Bernhard H. Neumann.

\par\noindent
Michel \textsc{Kervaire}, 1927--2007~;
fils de Heinz Hopf.

\par\noindent
Sergei Ivanovich \textsc{Adian}, 1931~; fils de Novikov.

\par\noindent
Michael Oser \textsc{Rabin}, 1931~; fils de Church.

\par\noindent
John Robert \textsc{Stallings}, 1935--2008~; fils de Fox.

\par\noindent
Sergei Petrovich \textsc{Novikov}, 1938~; fils de Postnikov.

\par\noindent
Friedhelm \textsc{Waldhausen}, 1938~;
fils de Hirzebruch.

\par\noindent
Andrew \textsc{Casson}, 1943~;
fils de C.T.C. Wall.

\par\noindent
Mikhael Leonidovich \textsc{Gromov}, 1943~;
fils de Rochlin.

\par\noindent
Gregori \textsc{Perelman}, 1966~;
fils d'Alexandrov et Burago.

\end{document}